\documentclass[11pt]{article}
\usepackage{amsmath, amssymb, amsfonts, amsthm, amscd, vmargin}
\usepackage[latin1]{inputenc}
\usepackage[all]{xy}
\usepackage{graphicx,color}
\usepackage{enumerate}
\usepackage{tikz}
\usetikzlibrary{shapes}

\setmarginsrb{2.4cm}{4.5cm}{2.4cm}{3.3cm}{0cm}{0cm}{0cm}{1.5cm}

\theoremstyle{plain}
\newtheorem{teo}{Theorem}
\newtheorem{lem}[teo]{Lemma}
\newtheorem{prop}[teo]{Proposition}
\newtheorem{cor}[teo]{Corollary}

\theoremstyle{definition}
\newtheorem{defi}{Definition}
\newtheorem{ex}{Example}

\newtheorem{rem}{Remark}
\newtheorem{rems}[rem]{Remarks}

\newcommand{\Hom}{\operatorname{Hom}}

\newcommand{\Ker}{\operatorname{Ker}}

\renewcommand{\dim}{\operatorname{dim}}
\newcommand{\codim}{\operatorname{codim}}

\renewcommand{\Im}{\operatorname{Im}}
\renewcommand{\Ker}{\operatorname{Ker}}

\newcommand{\Min}{\operatorname{Min}}

\newcommand{\Cbb}{{\mathbb C}}
\newcommand{\Qbb}{{\mathbb Q}}
\newcommand{\Zbb}{{\mathbb Z}}
\newcommand{\Rbb}{{\mathbb R}}
\newcommand{\Pbb}{{\mathbb P}}

\newcommand{\Fbb}{{\mathbb F}}

\newcommand{\lra}{\longrightarrow}

\begin{document}

\title{An approach of the Minimal Model Program for horospherical varieties via moment polytopes}

\author{Boris Pasquier\protect\footnote{Boris PASQUIER, Universit\'e Montpellier 2, CC 051, Place Eug\`ene Bataillon, 34095 Montpellier cedex, France. 
E-mail: boris.pasquier@math.univ-montp2.fr}}

\maketitle
\begin{abstract}
We describe the Minimal Model Program in the family of $\Qbb$-Gorenstein projective horospherical varieties, by studying a family of polytopes defined from the moment polytope of a Cartier divisor of the variety we begin with.
In particular, we generalize the results on MMP in toric varieties due to M.~Reid, and we complete the results on MMP in spherical varieties due to M.~Brion in the case of horospherical varieties.
\end{abstract}



~\\
\tableofcontents
	\section{Introduction}
	
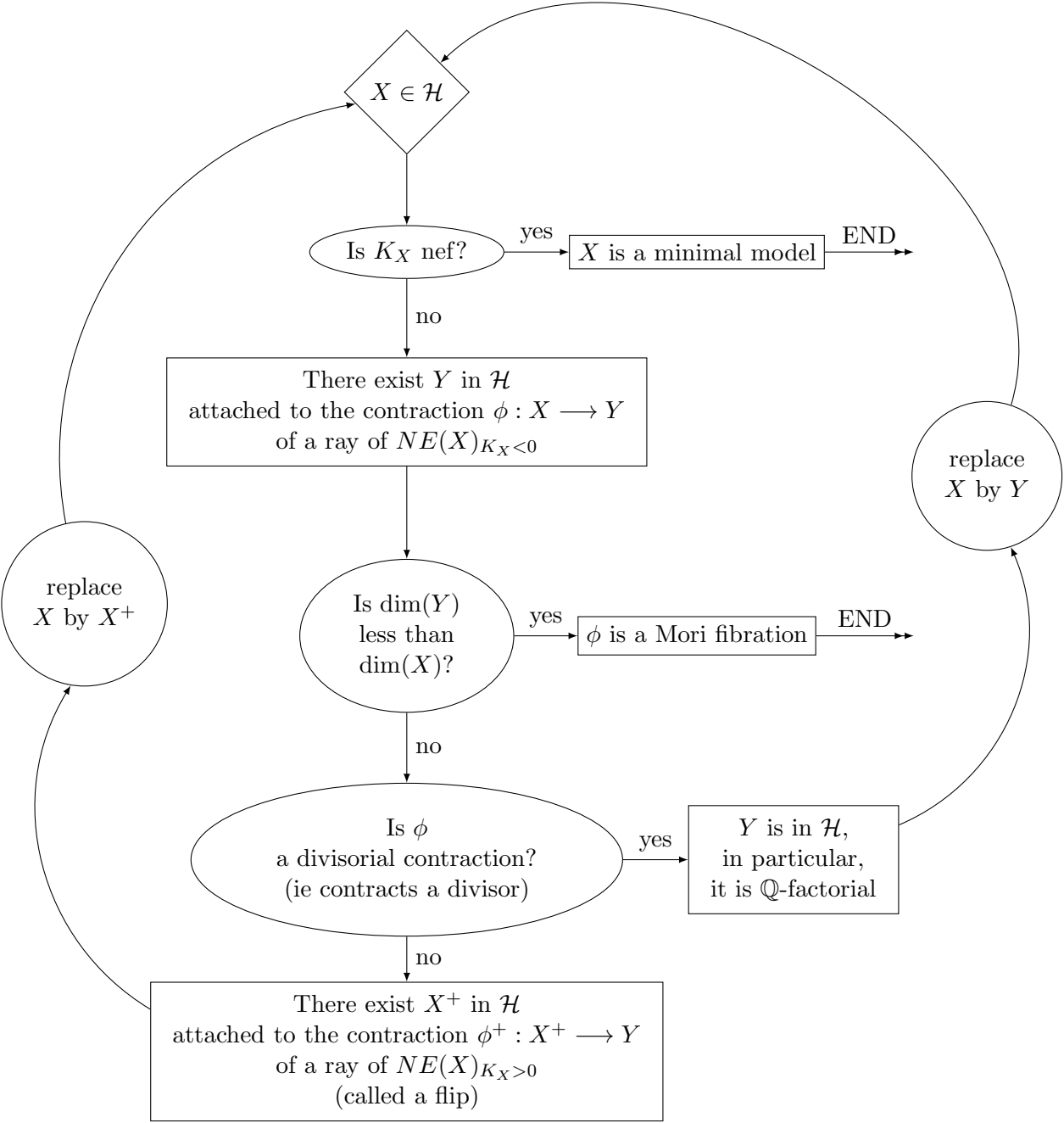
\begin{figure}
\caption{Original MMP for $\Qbb$-factorial varieties}
\label{fig:MMPoriginal}
\begin{tikzpicture}[scale=1]
\node [draw][diamond] (a) at (0,0) {
   $X\in\mathcal{H}$
};
\node [draw][ellipse] (b) at (0,-2.5) {
   Is $K_X$ nef?
};
\node [draw][rectangle] (c) at (4.5,-2.5) {
   $X$ is a minimal model
};
\node [rectangle] (d) at (8,-2.5) {};
\node [draw] [rectangle] (e) at (0,-5) {
\begin{tabular}{c}
There exist $Y$ in $\mathcal{H}$\\ attached to the contraction $\phi:X\longrightarrow Y$\\ of a ray of $NE(X)_{K_X<0}$
\end{tabular}
};
\node [draw] [ellipse] (f) at (0,-8.5) {
\begin{tabular}{c}
Is $\dim(Y)$\\ less than \\ $\dim(X)$?
\end{tabular}
};
\node [draw] [rectangle] (g) at (4.5,-8.5) {
$\phi$ is a Mori fibration
};
\node [rectangle] (h) at (8,-8.5) {};
\node [draw] [ellipse] (i) at (0,-12) {
\begin{tabular}{c}
Is $\phi$ \\ a divisorial contraction? \\(ie contracts a divisor)
\end{tabular}
};
\node [draw] [rectangle] (j) at (6,-12) {
\begin{tabular}{c}
$Y$ is in $\mathcal{H}$,\\ in particular,\\ it is $\Qbb$-factorial
\end{tabular}
};
\node [draw] [rectangle] (k) at (0,-15) {
\begin{tabular}{c}
There exist $X^+$ in $\mathcal{H}$\\ attached to the contraction $\phi^+:X^+\longrightarrow Y$\\ of a ray of $NE(X)_{K_X>0}$\\ (called a flip)
\end{tabular}
};
\node [draw][circle] (l) at (9,-6) {\begin{tabular}{c} replace \\ $X$ by $Y$\end{tabular}};
\node [draw][circle] (m) at (-5,-8) {\begin{tabular}{c} replace \\ $X$ by $X^+$\end{tabular}};

\draw[->,>=latex] (a) -- (b);
\draw[->,>=latex] (b) -- node[above] {yes} (c);
\draw[->>,>=latex] (c) --  node[above] {END} (d);
\draw[->,>=latex] (b) --  node[right] {no} (e);
\draw[->,>=latex] (e) -- (f);
\draw[->,>=latex] (f) --  node[above] {yes} (g);
\draw[->>,>=latex] (g) --  node[above] {END} (h);
\draw[->,>=latex] (f) --  node[right] {no} (i);
\draw[->,>=latex] (i) --  node[above] {yes} (j);
\draw[->,>=latex] (i) --  node[right] {no} (k);

\draw[->,>=latex] (j) to[bend right=45] (l);
\draw[->,>=latex] (l) to[bend right=75] (a);
\draw[->,>=latex] (k) to[bend left=45] (m);
\draw[->,>=latex] (m) to[bend left=45] (a);

\end{tikzpicture}
\end{figure}

The Minimal Model Program (MMP) takes an important place in birational algebraic geometry in order to get a birational classification of algebraic varieties. A lot of progress has been done in the last three decades. We come back here to an original version of the MMP, summarized by Figure~\ref{fig:MMPoriginal} where $\mathcal{H}$ denote a family of $\Qbb$-factorial varieties (see \cite{matsuki} to have a good overview of this theory).
For any $\Qbb$-Gorenstein variety $X$, we denote by $NE(X)$ the nef cone of curves of $X$, by $K_X$ a canonical divisor of $X$ and by $NE(X)_{K_X<0}$ (resp. $NE(X)_{K_X>0}$) the intersection of the nef cone with the open half-space of curves negative (resp. positive) along the divisor $K_X$.

When $\mathcal{H}$ is the family of $\Qbb$-factorial toric varieties, M.~Reid proved in 1983 that the MMP works \cite{Reid}. In particular, the MMP ends (there is no infinite series of flips). Note that, for this family, the cone $NE(X)$ is polyhedral generated by finitely many rays, and is very well-understood. Note also that, since toric varieties are rational varieties, a minimal model is then a point. Moreover, M.~Reid proved that the general fibers of the Mori fibrations, for $\Qbb$-factorial toric varieties, are weighted projective spaces.

When $\mathcal{H}$ is the family of $\Qbb$-factorial spherical $G$-varieties, for any connected reductive algebraic group $G$, M.~Brion proved in 1993 that the MMP works \cite{brionmori}. For this family, the cone $NE(X)$ is still polyhedral generated by finitely many rays and described in \cite{brionmori}. Note that spherical varieties are also rational varieties, so that minimal models are still points here. 
Nevertheless, it is very difficult to compute concretely $NE(X)$ and $K_X$, so that it makes difficult the application of the MMP to explicit examples of this family. That is why we reduce the study to horospherical varieties, for which a canonical divisor is well-known, and that is also why we present another approach that does not need the computation of $NE(X)$.
Moreover, the general fibers of Mori fibrations for $\Qbb$-factorial spherical $G$-varieties are not known. 

In this paper, we first consider the case where $\mathcal{H}$ is the family of $\Qbb$-Gorenstein projective horospherical $G$-varieties. The family of horospherical varieties is contained in the family of spherical varieties and contains toric varieties. Using a different approach than the one used by M.~Reid or M.~Brion, we obtain the main result of this paper. 

\begin{figure}
\caption{MMP for $\Qbb$-Gorenstein projective horospherical $G$-varieties}
\label{fig:MMPnew}
\begin{tikzpicture}[scale=1]

\node [draw][diamond] (a) at (0,0) {
   $X\in\mathcal{H}$
};
\node [draw][ellipse] (b) at (0,-2.5) {
   Is $K_X$ nef?
};
\node [draw][rectangle] (c) at (4.5,-2.5) {
   $X$ is a minimal model
};
\node [rectangle] (d) at (8,-2.5) {};
\node [draw] [rectangle] (e) at (0,-5) {
\begin{tabular}{c}
There exist $Y$ in $\mathcal{H}$\\ attached to the contraction $\phi:X\longrightarrow Y$\\ of a face of $NE(X)_{K_X<0}$
\end{tabular}
};
\node [draw] [ellipse] (f) at (0,-8.5) {
\begin{tabular}{c}
Is $\dim(Y)$\\ less than \\ $\dim(X)$?
\end{tabular}
};
\node [draw] [rectangle] (g) at (4.5,-8.5) {
$\phi$ is a Mori fibration
};
\node [rectangle] (h) at (8,-8.5) {};
\node [draw] [ellipse] (i) at (0,-12) {
\begin{tabular}{c}
Does $Y$ \\have 
$\Qbb$-Gorenstein\\ singularities? 
\end{tabular}
};
\node [draw] [rectangle] (j) at (6,-12) {
$\phi$ is a divisorial contraction 
};
\node [draw] [rectangle] (k) at (0,-15) {
\begin{tabular}{c}
There exist $X^+$ in $\mathcal{H}$\\ attached to the contraction $\phi^+:X^+\longrightarrow Y$\\ of a face of $NE(X)_{K_X>0}$
\end{tabular}
};
\node [draw][circle] (l) at (9,-6) {\begin{tabular}{c} replace \\ $X$ by $Y$\end{tabular}};
\node [draw][circle] (m) at (-5,-8) {\begin{tabular}{c} replace \\ $X$ by $X^+$\end{tabular}};

\draw[->,>=latex] (a) -- (b);
\draw[->,>=latex] (b) -- node[above] {yes} (c);
\draw[->>,>=latex] (c) --  node[above] {END} (d);
\draw[->,>=latex] (b) --  node[right] {no} (e);
\draw[->,>=latex] (e) -- (f);
\draw[->,>=latex] (f) --  node[above] {yes} (g);
\draw[->>,>=latex] (g) --  node[above] {END} (h);
\draw[->,>=latex] (f) --  node[right] {no} (i);
\draw[->,>=latex] (i) --  node[above] {yes} (j);
\draw[->,>=latex] (i) --  node[right] {no} (k);

\draw[->,>=latex] (j) to[bend right=45] (l);
\draw[->,>=latex] (l) to[bend right=75] (a);
\draw[->,>=latex] (k) to[bend left=45] (m);
\draw[->,>=latex] (m) to[bend left=45] (a);

\end{tikzpicture}
\end{figure}
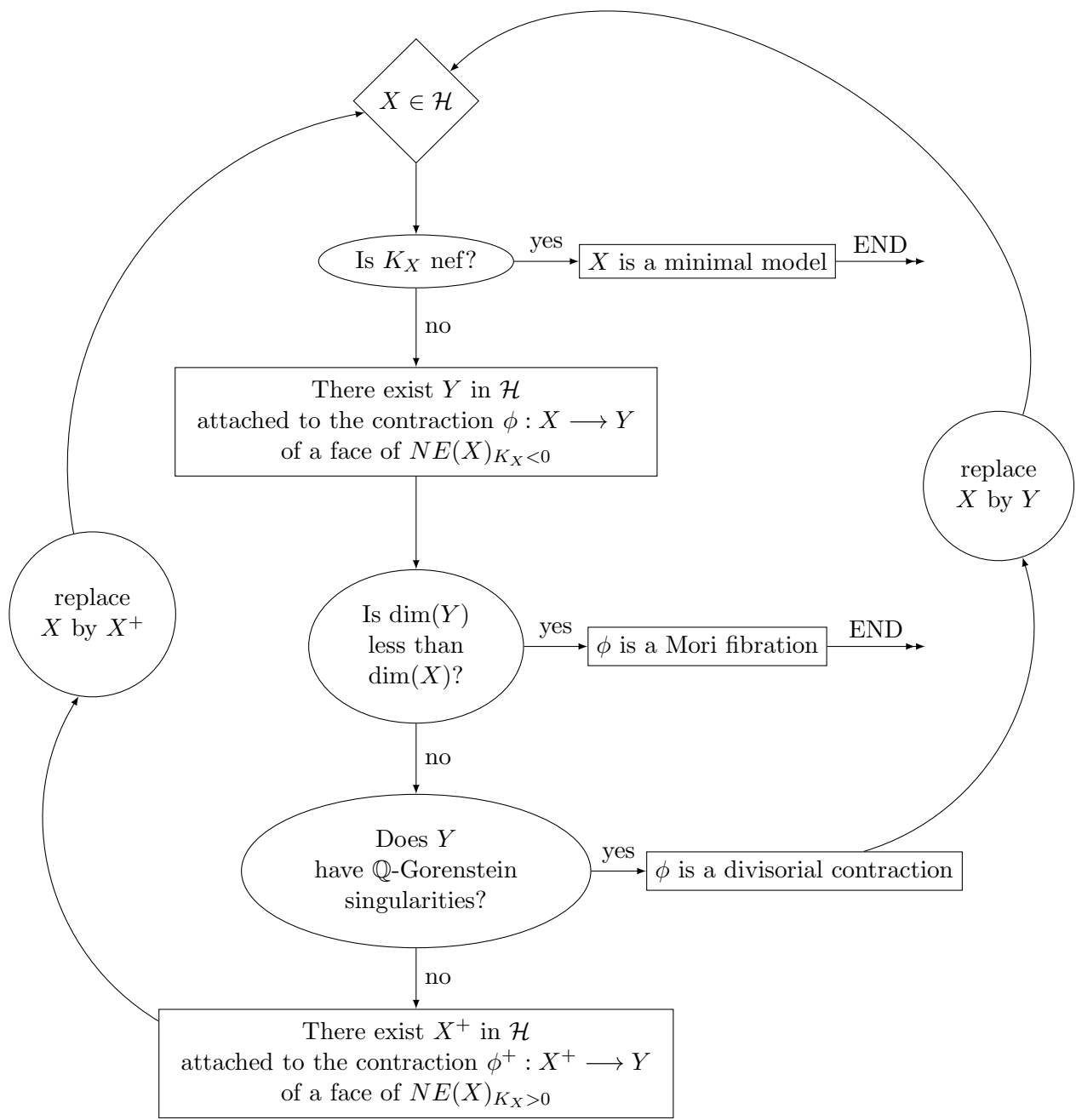

\begin{teo}\label{th:main0}
The MMP described by Figure~\ref{fig:MMPnew} works if $\mathcal{H}$ is the family of $\Qbb$-Gorenstein projective horospherical $G$-varieties for any connected reductive algebraic group $G$.
Moreover, for any $X$ in $\mathcal{H}$ and from any choice of an ample Cartier divisor of $X$,  we can concretely describe each step of this MMP until it ends.
\end{teo}

Remark that the main difference with the original MMP is that a divisorial contraction can give a flip,. It comes from the fact that the varieties are here not necessarily $\Qbb$-factorial.\\

We use here convex geometry, and in particular continuous transformations of moment polytopes. That is why we need to restrict to projective varieties. To illustrate Theorem~\ref{th:main0}, we give two first examples, both from the same toric variety $X$. Note that the result of the MMP applied to a variety $X$ is not unique: for the original MMP, it depends on the choices of a ray in the effective cone of all varieties appearing in the beginning of each loop of the program, and here it depends on the choice of a $\Qbb$-Cartier ample divisor of the variety $X$ (only in the beginning of the program). To choose a $\Qbb$-Cartier ample divisor of a projective toric variety is equivalent to choose its moment polytope (see Section~\ref{sec:momentpoly}).
We suppose here that the reader know the classification of toric varieties in terms of fans, see \cite{fulton} or \cite{oda} if it is not the case.

\begin{ex}\label{ex:toric1}
In $\Qbb^3$ consider the simple polytope $Q$ defined by the following inequations, where $(x,y,z)$ are the coordinates of a point in $\Qbb^3$.
$$\begin{array}{ccc}
z&\geq& -1\\
-x-y-2z&\geq &-5\\
2x-z&\geq& -3\\
-2x-z&\geq &-3\\
2y-z&\geq& -3\\
-2y-z&\geq& -3.
\end{array}$$

It is a pyramid whose summit is cut by a plane. 
 
The edges of the fan $\Fbb_X$ of the toric variety $X$ associated to $Q$ are $x_1:=(0,0,1)$, $x_2:=(-1,-1,-2)$, $x_3:=(2,0,-1)$, $x_4:=(-2,0,-1)$, $x_5:=(0,2,-1)$ and $x_6:=(0,-2,-1)$.
And the maximal cones of  $\Fbb_X$ are the cones respectively generated by $(x_1,x_3,x_5)$, $(x_1,x_3,x_6)$, $(x_1,x_4,x_5)$, $(x_1,x_4,x_6)$,  $(x_2,x_3,x_5)$, $(x_2,x_3,x_6)$, $(x_2,x_4,x_5)$, $(x_2,x_4,x_6)$.
(See Section~\ref{sec:momentpoly}, to have an explanation of the correspondence between moment polytopes and fans.)

Note that $X$ is $\Qbb$-factorial, because $Q$ is simple (and $\Fbb_X$ is simplicial).

We consider the family of polytopes $Q^\epsilon$ defines by the following inequations:
$$\begin{array}{ccc}
z&\geq& -1+\epsilon\\
-x-y-2z&\geq &-5+\epsilon\\
2x-z&\geq& -3+\epsilon\\
-2x-z&\geq &-3+\epsilon\\
2y-z&\geq& -3+\epsilon\\
-2y-z&\geq& -3+\epsilon.
\end{array}$$

Note that, for $\epsilon>0$ and small enough, it is the moment polytope of $D+\epsilon K_X$, where $D$ is the divisor of $X$ whose moment polytope is $Q$, and $K_X$ is the canonical divisor.

For $\epsilon\in[0,1[$, all polytopes $Q^\epsilon$ have the same structure so that they all correspond to the toric variety $X$. For any $\epsilon\in[1,2[$, the polytopes $Q^{\epsilon}$ are pyramids. They all correspond to the toric variety $Y$ whose fan $\Fbb_{Y}$ is described by the cones respectively generated by $(x_1,x_3,x_5)$, $(x_1,x_3,x_6)$, $(x_1,x_4,x_5)$, $(x_1,x_4,x_6)$,  $(x_3,x_4,x_5,x_6)$. Remark that the last cone is not simplicial, so that $Y$ is not $\Qbb$-factorial. The fact that $Y$ is $\Qbb$-Gorenstein comes from the fact that $x_3$, $x_4$, $x_5$ and $x_6$ are in a common plane. Note also that $x_2$ is not an edge of $\Fbb_{Y}$.

And for $\epsilon=2$, the polytope $Q^\epsilon$ is a point. Then the family $(Q^\epsilon)_{\epsilon\in[0,2]}$ reveals a divisorial contraction $\phi:X\longrightarrow Y$ and a Mori fibration from $Y$ to a point.

We illustrate this example in Figure~\ref{figure:exemple1}.

\begin{figure}[htbp]
\begin{center}
\input{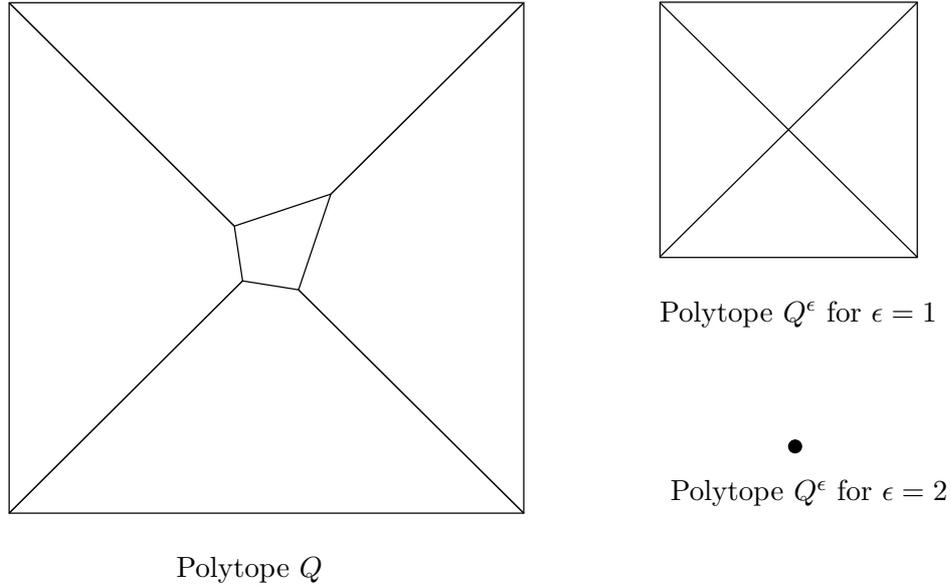}
\caption{Evolution of $Q^\epsilon$ in Example \ref{ex:toric1} (view from the top of the pyramid)}
\label{figure:exemple1}
\end{center}
\end{figure}
\end{ex}

Before to give the second example, we can notice that, in Example~\ref{ex:toric1}, the divisorial contraction $\phi$ goes from a $\Qbb$-factorial variety to a not $\Qbb$-factorial variety. It does not contradict the results of M.~Reid, because $\phi$ is not the contraction of a ray of $NE(X)$ but of a 2-dimensional face of $NE(X)$. That gives a reason of considering first $\Qbb$-Gorenstein varieties instead of $\Qbb$-factorial varieties.
Fortunately, what occurs in Example~\ref{ex:toric1} can be observed only in very particular cases (see Theorem~\ref{th:maingeneral}). Now, we consider a second example, with a more general divisor $D$.

\begin{ex}\label{ex:toric2}
 Let $Q^\epsilon$ be the polytopes defined by the following inequations,
$$\begin{array}{ccc}
z&\geq &-1+\epsilon\\
-x-y-2z&\geq& -5+\epsilon\\
2x-z&\geq& -4+\epsilon\\
-2x-z&\geq& -4+\epsilon\\
2y-z&\geq& -3+\epsilon\\
-2y-z&\geq& -3+\epsilon.
\end{array}$$
The simple polytope $Q^0$ is still a pyramid whose summit is cut by a plane. It is the moment polytope of a $\Qbb$-Cartier ample divisor of the toric variety $X$ of Example \ref{ex:toric1}.

For $\epsilon\in[0,\frac{1}{2}[$, all polytopes $Q^\epsilon$ have the same structure so that they all correspond to the toric variety $X$.

Now, for $\epsilon=\frac{1}{2}$, $Q^{\epsilon}$ corresponds to the fan whose maximal cones are the cones respectively generated by $(x_1,x_3,x_5)$, $(x_1,x_3,x_6)$, $(x_1,x_4,x_5)$, $(x_1,x_4,x_6)$,  $(x_2,x_4,x_5)$, $(x_2,x_4,x_6)$ and $(x_2,x_3,x_5,x_6)$. The associated toric variety $Y$ is not $\Qbb$-factorial. In fact, it is even not $\Qbb$-Gorenstein (because $x_2$, $x_3$, $x_5$ and $x_6$ are not in a common plane).

For $\epsilon\in]\frac{1}{2},\frac{3}{2}[$, the polytope $Q^{\epsilon}$ is the moment polytope of a divisor of the toric variety $X^+$ whose fan has maximal cones respectively generated by $(x_1,x_3,x_5)$, $(x_1,x_3,x_6)$, $(x_1,x_4,x_5)$, $(x_1,x_4,x_6)$,  $(x_2,x_4,x_5)$, $(x_2,x_4,x_6)$, $(x_2,x_5,x_6)$ and $(x_3,x_5,x_6)$. The variety $X^+$ is clearly $\Qbb$-factorial, and defines a flip $X\longrightarrow Y\longleftarrow X^+$.

Now, for $\epsilon\in[\frac{3}{2},2[$, $Q^{\epsilon}$ is a simple polytope with 6 vertices and the moment polytope of a divisor of the toric variety $Z$ whose fan has maximal cones respectively generated by  $(x_1,x_3,x_5)$, $(x_1,x_3,x_6)$, $(x_1,x_4,x_5)$, $(x_1,x_4,x_6)$,  $(x_3,x_5,x_6)$ and $(x_4,x_5,x_6)$. And we get a divisorial contraction from $X^+$ to $Z$. 

To finish, for $\epsilon=2$, $Q^{\epsilon}$ is the segment whose extremities are $(-\frac{1}{2},0,1)$ and $(\frac{1}{2},0,1)$, corresponding to the toric variety $\Pbb^1$. It gives a Mori fibration from $Z$ to the projective line.

We illustrate this example in Figure~\ref{figure:exemple2}.

\begin{figure}[htbp]
\begin{center}
\input{exemple2.pstex_t}
\caption{Evolution of $Q^\epsilon$ in Example \ref{ex:toric2} (view from the top of the pyramid)}
\label{figure:exemple2}
\end{center}
\end{figure}

Moreover, the general fiber of the Mori fibration is the toric variety associated to the polytope in $\Qbb^2$ defined by the inequalities $z\geq -1$, $2y-z\geq -3$ and $-2y-z\geq -3$. Its fan is isomorphic to the complete fan in $\Qbb^2$ whose edges are $(0,1)$, $(-2,1)$ and $(-2,-1)$. It is the weight projective plane $\Pbb(1,1,2)$. 
\end{ex}

Example~\ref{ex:toric2} arises a natural question: can we recover the original MMP for $\Qbb$-factorial projective horospherical varieties, by chosing a good divisor $D$? The answer is yes and given in the following result. 
\begin{teo}\label{th:maingeneral}
Suppose that $X$ is $\Qbb$-factorial.
Then, by taking a general divisor $D$ of $X$, we have that 
\begin{itemize}
\item at any loops of the MMP described in Figure~\ref{fig:MMPnew} (until it ends), the morphisms $\phi$ and $\phi^+$ are contractions of rays of $NE(X)$ and $NE(X^+)$;
\item the MMP described in Figure~\ref{fig:MMPnew} still works by replacing $\Qbb$-Gorenstein singularities with $\Qbb$-factorial ones everywhere;
\item the general fibers of Mori fibrations are projective $\Qbb$-factorial horospherical varieties with Picard number~1 (whose moment polytopes are simplexes intersecting all walls of a dominant chamber along facets).
\end{itemize}
\end{teo}

\begin{rem}
If $X$ is smooth and $D$ general, the general fibers of Mori fibrations are projective smooth horospherical varieties with Picard number~1. These varieties have been classified in \cite{2orbits}, in particular we get flag varieties and some two-orbit varieties.
Hence, the general fibers of Mori fibration are not only weighted projective spaces as in the toric case.
\end{rem}

The paper is organized as follows. 

In Section~\ref{sec:horopoly}, we recall the theory of horospherical varieties, we state the correspondence between polarized horospherical varieties and moment polytopes, we rewrite the existence criterion of an equivariant morphism between two horospherical varieties in terms of polytopes, and we describe the curves of horospherical varieties from a moment polytope. 

In Section~\ref{sec:poly}, we study particular (linear) one-parameter families of polytopes in $\Qbb^n$ and we define some equivalence relation in these families.

In Section~\ref{sec:MMP}, we prove Theorems~\ref{th:main0} and \ref{th:maingeneral}, by using the previous two sections. In particular, we construct a one-parameter families of moment polytopes whose equivalent classes describe all loops of the MMP given in Figure~\ref{fig:MMPnew}.

In Section~\ref{sec:examples}, we give five examples to illustrate what can happen in the MMP for horospherical (and not toric) varieties.

\section{Horospherical varieties and polytopes}\label{sec:horopoly}

\subsection{Notations and horospherical embedding theory}

We begin by recalling briefly the Luna-Vust theory of horospherical embeddings and by setting the notations used in the rest of the paper. For more details, the reader could have a look at \cite{Fanohoro}, or in more generalities at the Luna-Vust theory of spherical embeddings in \cite{knop}.

We fix a connected reductive algebraic group $G$. A closed subgroup $H$ of $G$ is said to be horospherical if it contains the unipotent radical $U$ of a Borel subgroup $B$ of $G$. It is equivalent to say that $G/H$ is a torus fibration over a flag variety $G/P$. The parabolic subgroup $P$ is the normalizer $N_G(H)$ of $H$ in $G$ and contains $B$. We fix a maximal torus $T$ of $B$. Then we denote by $S$ the set of simple roots of $(G,B,T)$. Also denote by $R$ the subset of $S$ of simple roots of $P$.
Let $X(T)$ (respectively $X(T)^+$) be the lattice of characters of $T$ (respectively the set of dominant characters). Similarly, we define $X(P)$ and $X(P)^+=X(P)\cap X(T)^+$. Note that $X(P)$ is generated by the fundamental weights $\varpi_\alpha$ with $\alpha\in S\backslash R$.
Let $M$ be the sublattice of $X(P)$ consisting of characters of $P$ vanishing on $H$. The rank of $M$ is called the rank of $G/H$ and denoted by $n$. Let $N:=\Hom_\Zbb(M,\Zbb)$. For any free lattice $\mathbb{L}$, we denote by $\mathbb{L}_\Qbb$ the $\Qbb$-vector space $\mathbb{L}\otimes_\Zbb\Qbb$.  For any simple root $\alpha\in S\backslash R$, the restriction of the coroot $\alpha^\vee$ to $M$ is a point of $N$, which we denote by $\alpha^\vee_M$. For any $\alpha\in S\backslash R$, we also denote by $W_{\alpha,P}$ the hyperplane defined by $\{m\in X(P)_\Qbb\,\mid\,\langle m,\alpha^\vee\rangle=0\}$ (note that it corresponds to a wall of the dominant chamber of characters of $P$).

A $G/H$-embedding is couple $(X,x)$, where $X$ is a normal algebraic $G$-variety and $x$ a point of $X$ such that $G\cdot x$ is open in $X$ and isomorphic to $G/H$. By abuse of notation, we will forget the point $x$.

The $G/H$-embeddings are classified by colored fans in $N_\Qbb$. A colored fan is a set of colored cones in $N_\Qbb$ stable by taking colored subcones and without overlap (see for example \cite{Fanohoro} for the complete definition of a colored fan). A colored cone is a couple $(\mathcal{C},\mathcal{F})$, where $\mathcal{F}$ is a subset of $S\backslash R$ such that $\alpha^\vee_M\neq 0$ for all $\alpha\in\mathcal{F}$, and $\mathcal{C}$ is a strictly convex cone generated by the $\alpha^\vee_M$ with $\alpha\in\mathcal{F}$ and a finite set of points in $N$. 

If $G=(\Cbb^*)^n$ and $H=\{1\}$, we get the well-known classification of toric varieties.

The colored fans of complete $G/H$-embeddings are the complete ones (ie such that $N_\Qbb$ is the union of the colored cones of the colored fan).

If $X$ is a $G/H$-embedding, we denote by $\mathbb{F}_X$ the colored fan of $X$ in $N_\Qbb$ and we denote by $\mathcal{F}_X$ the set $\cup_{(\mathcal{C},\mathcal{F})\in\mathbb{F}_X}\mathcal{F}\subset S\backslash R$ of colors of $X$.

Moreover, the set of $G$-orbits of a $G/H$-embedding $X$ is naturally in bijection with the set of colored cones of $\mathbb{F}_X$ (reversing usual orders). In particular, the $G$-stable irreducible divisors correspond to the colored edges of $\mathbb{F}_X$ of the form $(\mathcal{C},\varnothing)$. We denote them by $X_1,\dots,X_m$, and for any $i\in\{1,\dots,m\}$, we denote by $x_i$ the primitive element in $N$ of the corresponding edge.

\subsection{Correspondence between projective horospherical varieties and polytopes}\label{sec:momentpoly}

In this section, we list results coming directly from the characterization of Cartier, $\Qbb$-Cartier and ample divisors of horospherical varieties due to M.~Brion in the more general case of spherical varieties (\cite{briondiv}). And we classify projective $G/H$-embeddings in terms of $G/H$-polytopes (defined below in Definition~\ref{def:G/H-eq}).

First, describe the $B$-stable irreducible divisors of a $G/H$-embedding $X$. We already defined the $G$-stable divisors $X_1,\dots,X_m$. The other ones are the closures in $X$ of $B$-stable irreducible divisors of $G/H$, which are the inverse images by the torus fibration $G/H\longrightarrow G/P$ of the 
Schubert divisors of the flag variety $G/P$. We denote these divisors of $X$ by $D_\alpha$ for any $\alpha\in S\backslash R$.

In all the paper, a divisor of a horospherical variety is always supposed to be $B$-stable.

Let $X$ be a projective $G/H$-embedding and let $D=\sum_{i=1}^m a_i X_i +\sum_{\alpha\in S\backslash R} a_\alpha D_\alpha$ be an ample Cartier divisor of $X$. Then the set $$\tilde{Q}_{X,D}=\tilde{Q}_D:=\{m\in M_\Qbb\,\mid\,\langle m,x_i\rangle\geq -a_i,\,\forall i\in\{1,\dots,m\}\mbox{ and }\langle m,\alpha^\vee_M\rangle\geq -a_\alpha,\,\forall \alpha\in\mathcal{D}_X\}$$ is a lattice polytope in $M_\Qbb$ (ie with vertices in $M$) of maximal dimension. If $D$ is only $\Qbb$-Cartier, then $\tilde{Q}_D$ is a polytope in $M_\Qbb$ of maximal dimension. The polytope $Q_{X,D}=Q_D:=v^0+\tilde{Q}_D$ in $X(P)_\Qbb$, where $v^0=\sum_{\alpha\in S\backslash R}a_\alpha\varpi_\alpha$,  is called the moment polytope of the polarized variety $(X,D)$ (or of $D$ in $X$). We call $\tilde{Q}_D$ the pseudo-moment polytope of $(X,D)$. 

Moreover, the moment polytope of an ample $\Qbb$-Cartier divisor $D$ in a $G/H$-embedding $X$ is always contained in the dominant chamber $X(P)^+$, and $Q_D\cap W_{\alpha,P}\neq\emptyset$ if and only if $\alpha\in\mathcal{F}_X$. We note also that, since any $\alpha\in\mathcal{F}_X$ satisfies $\alpha_M^\vee\neq 0$, a moment polytope is contained in no wall $W_{\alpha,P}$, with $\alpha\in S\backslash R$.

An important tool of the paper is the fact that the colored fan of $X$ is reconstructible from the moment polytope $Q_D$ of $(X,D)$ for any ample $\Qbb$-Cartier ($B$-stable) divisor $D$. 
Indeed, any maximal colored cone of the complete colored fan $\mathbb{F}_X$ can be defined from a vertex of $Q_D$ as follows. Let $v$ be a vertex of $Q_D$. We define $\mathcal{C}_v$ to be the cone of $N_\Qbb$ generated by inward-pointing normal vectors of the facets of $Q_D$ containing $v$. And we set $\mathcal{F}_v=\{\alpha\in S\backslash R\,\mid\,v\in W_{\alpha,P}\}$. Then $(\mathcal{C}_v,\mathcal{F}_v)$ is a maximal colored cone of $\mathbb{F}_X$.

Moreover, the divisor $D$ can also be computed from the pair $(Q,\tilde{Q})$. Indeed, the coefficients $a_\alpha$ with $\alpha\in S\backslash R$ are given by the translation vector in $X(P)^+$ that maps $\tilde{Q}$ to $Q$. And, for any $i\in\{1,\dots,m\}$, the coefficient $a_i$ is given by $\langle v_i,x_i\rangle$ for any element $v_i\in M_\Qbb$ in the facet of $\tilde{Q}$ for which $x_i$ is an inward-pointing normal vector.

In order to classify projective $G/H$-embeddings in terms of polytopes, we give the following definition.

\begin{defi}\label{def:G/H-eq}
Let $Q$ be a polytope in $X(P)^+_\Qbb$. We say that $Q$ is a $G/H$-polytope, if its direction is $M_\Qbb$ and if it is contained in no wall $W_{\alpha,P}$ with $\alpha\in S\backslash R$.

Let $Q$ and $Q'$ be two $G/H$-polytopes in $X(P)^+_\Qbb$. Consider any polytopes $\tilde{Q}$ and $\tilde{Q'}$ in $M_\Qbb$ obtained by translations of $Q$ and $Q'$ respectively.   We say that $Q$ and $Q'$ are equivalent $G/H$-polytopes if the following conditions are satisfied.
\begin{enumerate}
\item There exist an integer $j$, $2j$ affine half-spaces $\mathcal{H}_1^+,\dots,\mathcal{H}_j^+$ and $\mathcal{H'}_1^+,\dots,\mathcal{H'}_j^+$ of $M_\Qbb$ (respectively delimited by the affine hyperplanes $\mathcal{H}_1,\dots,\mathcal{H}_j$ and $\mathcal{H'}_1,\dots,\mathcal{H'}_j$) such that $\tilde{Q}$ is the intersection of the $\mathcal{H}_i^+$, $\tilde{Q'}$ is the intersection of the $\mathcal{H'}_i^+$, and for all $i\in\{1,\dots,j\}$, $\mathcal{H}_i^+$ is the image of $\mathcal{H'}_i^+$ by a translation.
\item With notations of the previous item, for all subset $J$ of $\{1,\dots,j\}$, the intersections $\cap_{i\in J}\mathcal{H}_i\cap Q$ and $\cap_{i\in J}\mathcal{H'}_i\cap Q'$ have the same dimension.
\item $Q$ and $Q'$ intersect exactly the same walls $W_{\alpha,P}$ of $X(P)^+$ (with $\alpha\in S\backslash R$).
\end{enumerate}
\end{defi}

Remark that this definition does not depend on the choice of $\tilde{Q}$ and $\tilde{Q'}$.
As a corollary of what we say just above, we obtain the following result.

\begin{prop}\label{prop:classprojembpoly}
The correspondence between moment polytopes and colored fans gives a bijection between the set of classes of $G/H$-polytopes and (isomorphic classes of) projective $G/H$-embeddings.

Moreover, the set of $G$-orbits of a projective $G/H$-embedding is in bijection with the set of faces of one of its moment polytope (preserving the respective orders). 
\end{prop}

In section~\ref{sec:generalfiber}, we need precise description of the $G/H$-embedding associated to a $G/H$-polytope.
For any dominant weight $\chi$, we denote by $V(\chi)$ the irreducible $G$-module of highest weight $\chi$, and we fix a highest weight vector $v_\chi$ in $V(\chi)$. The Borel subgroup of $G$  opposite to $B$ is denoted by $B^-$. 

\begin{prop}\label{prop:veryample}
Suppose $D$ is Cartier and very ample, then $X$ is isomorphic to the closure of $G\cdot [\sum_{\chi\in (v^0+M)\cap Q}v_\chi]$ in $\Pbb(\oplus_{\chi\in (v^0+M)\cap Q}V(\chi))$.
\end{prop}

\begin{rem}\label{rem:veryample}
If $D$ is Cartier and only ample, then $(n-1)D$ is very ample (\cite[Theorem 0.3]{these}). In particular, the assumptions Cartier and very ample, instead of $\Qbb$-Cartier and ample, are not really restrictive. Indeed, if $D$ is $\Qbb$-Cartier and ample, then there exists a non-zero integer $p$ such that $pD$ is Cartier and very ample, and $X$ is isomorphic to the closure of $G\cdot [\sum_{\chi\in (pv^0+M)\cap pQ)}v_\chi]$ in $\Pbb(\oplus_{\chi\in (pv^0+M)\cap pQ}V(\chi))$.
\end{rem}

\begin{proof}
First, since $\tilde{Q}$ is lattice polytope of maximal dimension, the $G$-orbit $G\cdot [\sum_{\chi\in (v^0+M)\cap Q}v_\chi]$ is isomorphic to $G/H$. By \cite[Lemma~5.1]{these}, the closure $X'$ of this $G$-orbit is normal, and then a $G/H$-embedding. Now, if $\chi_0$ is a vertex of $Q$, the intersection of the closure of the $G$-orbit with the affine space $\cap\Pbb(\oplus_{\chi\in (pv^0+M)\cap pQ}V(\chi))_{v_{\chi_0}^*\neq 0}$ is the $B^-$-stable affine variety whose structure ring is the $B^-$-module generated by the rational functions $\frac{v_\chi^*}{v_{\chi_0}^*}$ of $\Pbb(\oplus_{\chi\in (v^0+M)\cap Q}V(\chi))$,  with $\chi_0\neq\chi\in (v^0+M)\cap Q$. It is also the $B^-$-stable affine variety constructed from the maximal colored cone of $\mathbb{F}_X$ associated to the vertex $\chi_0$ of $Q$. Hence, the colored fan of $X$ is the same as the colored fan of $X'$ and $X$ is isomorphic to $X'$.
\end{proof}

\subsection{$G$-equivariant morphisms and polytopes}\label{sec:morph}

The existence of $G$-equivariant morphisms between horospherical varieties can be characterized in terms of colored fans \cite{knop}.
In this section, we rewrite this characterization in terms of moment polytopes.

Let $(X,D)$ be a polarized $G/H$-embedding, let $(X',D')$ be a polarized $G/H'$-embedding and denote by $Q$ and $Q'$ the corresponding moment polytopes respectively. (We denote all data corresponding to $G/H'$ and $X'$ with primes: $H'$, $P'$, $R'$, $M'$, $N'$,...). We denote by $\tilde{Q}$ and $\tilde{Q'}$ the pseudo-moment polytopes.

A first necessary condition to the existence of a dominant $G$-equivariant morphism from $X$ to $X'$, is that there exists a projection $\pi$ from $G/H$ to $G/H'$. In particular $H'\supset H$, $P'\supset P$, $R'\supset R$. The projection $\pi$ induces an injective morphism $\pi^*$ from $M'$ to $M$ and a surjective morphism $\pi_*$ from $N$ to $N'$. We suppose that this necessary condition is satisfied in all the rest of the section and we identify $M'$ with $\pi_*(M')$.

Now, we define an application $\psi$ from the set of facets of $\tilde{Q}$ to the set of faces of $\tilde{Q'}$.
First, note a general fact on polytopes: if $\mathcal{P}$ is a polytope in $\Qbb^r$, then for any affine half-space $\mathcal{H}^+$ delimited by the affine hyperplane $\mathcal{H}$ in $\Qbb^r$, there exists a unique face $F$ of $\mathcal{P}$ such that there exists $x\in\Qbb^r$ such that $F$ is defined by $x+\mathcal{H}$ (ie $F=\mathcal{P}\cap (x+\mathcal{H})$ and $\mathcal{P}\subset v+ \mathcal{H}^+$).
Then, for any facet $F$ of $\tilde{Q}$, let  $\mathcal{H}^+$ be the affine half-space in $M_\Qbb$ containing $F$. If $\mathcal{H}^+\cap M'_\Qbb\neq M'_\Qbb$, it is an affine half-space in $M'_\Qbb$ and, applying the fact above to $\mathcal{P}=\tilde{Q'}$, it gives a unique face $F'$ of $\tilde{Q'}$. We set $\psi(F)=F'$. And if $\mathcal{H}^+\cap M'_\Qbb=M'_\Qbb$,  we  set $\psi(F)=\tilde{Q'}$.

\begin{prop}\label{prop:morph}
Under the above conditions, there exists a dominant $G$-equivariant morphism from $X$ to $X'$, if and only if 
\begin{enumerate}
\item for any subset $\mathcal{G}$ of facets of $\tilde{Q}$, $\cap_{F\in\mathcal{G}}F\neq\emptyset$ implies  $\cap_{F\in\mathcal{G}}\psi(F)\neq\emptyset$, and
\item for any $\alpha\in S\backslash R$ such that $Q\cap W_{\alpha,P}\neq\emptyset$ , we have $Q'\cap W_{\alpha,P}\neq\emptyset$. 
\end{enumerate}
\end{prop}

Remark that, in Proposition~\ref{prop:morph}, we can replace $\tilde{Q}$ by $Q$ (by extending the definition of $\psi$).

\begin{proof}
For the second condition, remark that if $\alpha\in S\backslash R$ is in $R'$, then $Q'$ is contained in $W_{\alpha,P}$. And, if $\alpha\in S\backslash R'$, $Q'\cap W_{\alpha,P}\neq\emptyset$ is equivalent to $Q'\cap W_{\alpha,P'}\neq\emptyset$  because $Q'\subset X(P')$ and $W_{\alpha,P'}=X(P')\cap W_{\alpha,P}$. Recall also that $Q\cap W_{\alpha,P}\neq\emptyset$ is equivalent to $\alpha\in\mathcal{F}_{X}$ (and $Q'\cap W_{\alpha,P'}\neq\emptyset$ is equivalent to $\alpha\in\mathcal{F}_{X'}$).

By \cite{knop}, there exists a dominant $G$-equivariant morphism from $X$ to $X'$ if and only if, for all colored cone $(\mathcal{C},\mathcal{F})$ of $X$, there exists a colored cone $(\mathcal{C'},\mathcal{F'})$ of $X'$ such that $\pi_*(\mathcal{C})\subset\mathcal{C'}$ and any element $\alpha\in\mathcal{F}$, either $\alpha\in R'$ or $\alpha\in\mathcal{F'}$.

In our case, since $X$ and $X'$ are complete, we can rewrite this characterization as follows. Denote by $y_i$  with $i\in\{1,\dots,k\}$ (resp. $y_i'$ with $i\in\{1,\dots,k'\}$) the primitive elements of the egdes of the colored fan of $X$ (resp. $X'$). For all $j\in\{1,\dots,k\}$, let $J'_j$ be the minimal subset of $\{1,\dots,k'\}$ such that $\pi_*(y_j)$ is in the cone $\mathcal{C'}_j$ generated by the $y'_i$ with $i\in J'_j$ (it is never empty because $X'$ is complete). We prove now that there exists a dominant $G$-equivariant morphism from $X$ to $X'$ if and only if, for all (colored) cone $\mathcal{C}$ of $X$, generated by the $y_j$ with $j\in J_\mathcal{C}$, the cone generated by the $y'_i$ with $i\in\cup_{j\in J_\mathcal{C}}J'_j$ is contained in a (colored) cone of $X'$; and $\mathcal{F}_X$ is contained in $R'\cup\mathcal{F}_{X'}$.

Indeed, let $\mathcal{C}$ be a cone of $X$, generated by the $y_j$ with $j\in J_\mathcal{C}$. If there exists a cone $\mathcal{C'}$ of $X'$ such that $\pi_*\mathcal{C}\subset\mathcal{C'}$, then $\mathcal{C'}$ contains all $y_j$ with $j\in J_\mathcal{C}$. Then by minimality of the $J'_j$, it contains also cones $\mathcal{C'}_j$ with $j\in J_\mathcal{C}$ and also the cone generated by the $y'_i$ with $i\in\cup_{j\in J_\mathcal{C}}J'_j$.
Conversely, if the cone generated by the $y'_i$ with $i\in\cup_{j\in J_\mathcal{C}}J'_j$ is contained in a cone $\mathcal{C'}$ of $X'$, then it is obvious that $\pi_*(\mathcal{C})\subset\mathcal{C'}$
And the condition on colors is the same in both cases.

Now, the proposition comes from the bijective correspondence between colored cones of $X$ (resp.  $X'$) and faces of $Q$ (resp. $Q'$), the first paragraph of the proof, and the following fact: the intersection of some facets of $Q$ is not empty if and only if the cone generated by the inward-pointing normal vectors corresponding to these facets is included in a cone of $X$ (and this latter cone corresponds to the face defined as the intersection of these facets).
\end{proof}

\begin{cor}
Suppose there exists a dominant $G$-equivariant morphism $\phi$ from $X$ to $X'$. Let $\mathcal{O}$ be the $G$-orbit in $X$ associated to a face $\cap_{F\in\mathcal{G}}F$. Then $\phi(\mathcal{O})$ is the $G$-orbit in $X'$ associated to a face $\cap_{F\in\mathcal{G}}\psi(F)$.
\end{cor}

\begin{proof}
Rewrite, in terms of polytopes, the fact that if $\mathcal{O}$ is the $G$-orbit in $X$ associated to a colored cone $(\mathcal{C},\mathcal{F})$ then $\phi(\mathcal{O})$ is the $G$-orbit in $X'$ associated to the minimal colored cone $(\mathcal{C},\mathcal{F})$ such that $\pi_*(\mathcal{C})\subset\mathcal{C'}$.
\end{proof}

\subsection{Curves in horospherical varieties}\label{sec:curves}

We begin this section by collecting, in the following Theorem, some results on curves in spherical varieties due to M.~Brion \cite{brionmori}. 
We denote by $N_1(X)$ the group of numerical classes of 1-cycles of the variety $X$. Recall that $NE(X)$ is the convex cone in $N_1(X)$ generated by effective 1-cycles.
We denote by $N^1(X)$ the group of numerical classes of $\Qbb$-Cartier divisors of $X$, it is the dual of $N_1(X)$.

\begin{teo}[M.~Brion]
Let $X$ be a complete spherical variety. Let $Y$ and $Z$ be two distinct closed $G$-orbits in $X$. Note that $Y$ and $Z$ are flag varieties, so that they have exactly one point fixed by $B$.
\begin{enumerate}
\item There exists a $B$-stable curve $C$ containing the $B$-fixed points of $Y$ and $Z$ if and only if the colored cones in the colored fan $\Fbb_X$ corresponding to the two closed $G$-orbits $Y$ and $Z$ intersect along a one-codimensional colored cone $\mu$). And, in that case, $C$ is unique, isomorphic to $\Pbb^1$, we denote it by $C_\mu$.
\item The space $N_1(X)_\Qbb$ is generated by the classes of the curves $C_\mu$ where $\mu$ is any one-codimensional colored cone in $\Fbb_X$, and some classes $[C_{D,Y}]$ where $Y$ is any closed $G$-orbit of $X$ and $D$ is any irreducible $B$-stable divisor of $X$ not containing $Y$. 
\item Any classes $[C_{D,Y}]$ that is not in a ray generated by the class of a curve $C_\mu$ is represented by a (not unique in general) $B$-stable curve $C_{D,Y}$ admitting the $B$-fixed point of $Y$ as unique $B$-fixed point. 
\item The cone $NE(X)$ is generated by the classes $[C_\mu]$  and $[C_{D,Y}]$.
\end{enumerate}
\end{teo}

Now, in the particular case of horospherical varieties, we complete this collection of results by the following proposition that describes an explicit curve in all classes of the form $[C_{D,Y}]$.
Note that a homogeneous $G$-space has a $B$-fixed point if and only if it is complete (projective), so that the $B$-fixed points of a complete $G$-variety are the (unique) $B$-fixed points of its closed $G$-orbits.
Denote by $s_\alpha$ the simple reflection associated to a simple root $\alpha$, and denote by $w_0$ the longest element of the Weyl group of $(G,T)$. 

\begin{prop} \label{prop:curves}
\begin{enumerate}
\item Let $X$ be a complete horospherical variety.
If a $B$-stable curve $C$ in $X$ contains a unique $B$-fixed point $y$, then $C$ is contained in the closed $G$-orbit $Y:=G\cdot y$. In particular, it is a Schubert subvariety of $Y$. 

\item For any closed $G$-orbit $Y$ and any divisor $D_\alpha$ that does not contain $Y$, the class of $[C_{D_\alpha,Y}]$ is represented by the Schubert variety of $Y$ given by the simple root $\alpha$ (ie $\overline{Bs_\alpha\cdot y}$, where $y$ is the $B$-fixed point of $Y$). We denote it by $C_{\alpha,Y}$.
\end{enumerate}
\end{prop}

\begin{proof}
\begin{enumerate}
\item Let $C$ be a $B$-stable curve in $X$. Suppose that $C$ is not contained in $Y$. By replacing $X$ by the closure of the biggest $G$-orbit of $X$ that intersects $C$, we can assume that $C$ intersect the open $G$-orbit $G/H$ that is at least of rank one. And then the intersection $C\cap G/H$ is an open set of $C$. Recall that $G/H$ is a $G$-equivariant torus fibration over the flag variety $G/P$. 
Then $C\cap G/H$ is the fiber $P/H$ of this fibration over $P/P$. Indeed, the $B$-orbits of $G/H$ are the inverse image of the Schubert cells of $G/P$. The (unique) smallest one is $P/H$, but by hypothesis, $P/H$ has positive dimension. Then, since $C\cap G/H$ is $B$-stable, it is $P/H$. Moreover $P/H$ has to be one-dimensional.

Prove now that $C$ has two $B$-fixed points. By the previous paragraph, $P/H$ is one-dimensional, ie $X$ is of rank one. Then, there exists a $\Pbb^1$-bundle $\tilde{X}$ over $G/P$ and a $G$-equivariant birational morphism $\phi:\tilde{X}\lra X$ ($\tilde{X}$ is the toroidal variety over $X$, see \cite[Example 1.13 (2)]{these}). The closure of the $P$-orbit $P/H$ in $\tilde{X}$ is the fiber $\Pbb^1$ over $P/P$. Then it has exactly two $B$-fixed points corresponding to the two $\Cbb^*$-fixed points of the toric variety $\Pbb^1$. Moreover, the two-closed $G$-orbits of $\tilde{X}$ are send, by $\phi$, respectively to the two-closed $G$-orbits of $X$ (here we use that $X$ is horospherical of rank one). Then $C=\phi(\tilde{C})$ has also two $B$-fixed points (and is isomorphic to $\Pbb^1$).

\item Let  $Y$ be a $G$-closed orbit of $X$ and let $\alpha\in S\backslash R$ such that $D_\alpha$ does not contain $Y$. Let $C_{\alpha,Y}:=\overline{Bs_\alpha\cdot y}$, where $y$ is the $B$-fixed point of $Y$. 
Let $D=\sum_{i=1}^mb_iX_i+\sum_{\beta\in S\backslash R}b_\beta D_\beta$ be a Cartier divisor of $X$. We want to compute $D\cdot C_{\alpha,Y}$.
Since $D$ is Cartier, there exists an eigenvector $f_Y(D)$ of $B$ in the set of rational function on $X$ such that the support of $D-\rm{div}(f_Y(D))$ is in the union of irreducible $B$-stable divisors that do not contain $Y$. Denote by $\chi_Y(D)$ the weight of the eigenvector $f_Y(D)$. But for any divisor $X_i$ with $i\in\{1,\dots,m\}$ such that $Y\not\subset X_i$, we clearly have $X_i\cdot C_{\alpha,Y}=0$ because $Y$ and $X_i$ are disjoint. Hence $$D\cdot C_{\alpha,Y}=(D-\rm{div}(f_Y(D)))\cdot C_{\alpha,Y}=\sum_{\beta,\,Y\not\subset D_\beta}(b_\beta-\langle\chi_Y(D),\beta^\vee_M\rangle)D_\beta\cdot C_{\alpha,Y}.$$
We conclude, by the following Lemma \ref{lem:intersection} and the formula given in \cite[3.2]{brionmori}.
\end{enumerate}
\end{proof}

\begin{lem}\label{lem:intersection}
Let $Y$ be a $G$-closed orbit of $X$ and let $\alpha,\,\beta\in S\backslash R$ such that $D_\alpha$ and $D_\beta$ do not contain $Y$. 
Then $D_\beta\cdot C_{\alpha,Y}=\delta_{\alpha\beta}$, where $\delta$ is the Kronecker delta.
\end{lem}

\begin{proof}
Denote by $y$ the $B$-fixed point of $Y$. Remark that $y$ is also fixed by $P$.
By \cite[Lemma 2.8]{these}, there exists a unique $G$-equivariant morphism $\phi$ from the open $G$-stable set $X_Y:=\{x\in X\mid\overline{G\cdot x}\supset Y\}$ to the closed $G$-orbit $Y$. Note that  $y=\phi(H/H)$ because $y$ is the only point of $Y$ fixed by $U\subset H$ and, similarly, $\phi$ is the identity on $Y$.

Recall that $D_\beta$ is the closure in $X$ of the $B$-orbit $Bs_\beta w_0P/H$. Then $D_\beta\cap Y$ is contained in (and then equals) the closure of $\phi(Bs_\beta w_0P/H)$ in $Y$, which is the Schubert variety $\overline{Bs_\beta w_0\cdot y}$ of $Y$. Moreover, since $D_\beta$ does not contain $Y$, $D_\beta\cap Y$ is a divisor of $Y$. Hence $D_\beta\cdot C_{\alpha,Y}=(D_\beta\cap Y)\cdot C_{\alpha,Y}=\delta_{\alpha\beta}$. 
\end{proof}

With the correspondence between colored fans and moment polytopes, if $X$ is a horospherical variety and $D$ is a $\Qbb$-Cartier divisor, we denote by $C_\mu$ for any edge $\mu$ of $Q_D$, and by $C_{\alpha,v}$ for any $\alpha\in S\backslash R$ and any vertex $v$ not contained in $W_{\alpha,P}$, the curves defined above in terms of the colored fans (one-codimensional colored cones correspond to edges of moment polytopes, and closed $G$-orbits correspond to maximal colored cones and then to vertices of moment polytopes).

The following results is a direct consequence of the  formula of \cite[section 3.2]{brionmori}, already used in the proof of Proposition~\ref{prop:curves}.

\begin{prop}\label{prop:intersection}
Let $X$ be a horospherical variety and $D$ an ample $\Qbb$-Cartier divisor. 

Then, for any edge $\mu$ of the moment polytope $Q_D$, the intersection number $D.C_\mu$ is the integral length of $\mu$, ie the length of $\mu$ divided by the length of the primitive element in the direction of $\mu$. 

And for any $\alpha\in S\backslash R$ and for any vertex $v$ of $Q_D$ not in the wall $W_{\alpha,P}$, we have $D.C_{\alpha,v}=\langle v,\alpha^\vee\rangle$.
\end{prop}

\section{One-parameter families of polytopes}\label{sec:poly}

In this section, we study particular one-parameter families of polytopes. This section can be read independently from the rest of the paper. Corollary~\ref{cor:mainpoly} is an essential tool in the proofs of Theorems~\ref{th:main0} and \ref{th:maingeneral}.

\subsection{A first one-parameter family of polytopes: definitions and results}

First, we fix notations. Let $n$ and $m$ be two positive integers. Consider three matrices $A$, $B$ and $C$ respectively in $M_{m\times 1}(\Qbb)$, $M_{m\times 1}(\Qbb)$ and $M_{m\times n}(\Qbb)$. Then we define a first family of polyhedrons indexed by $\epsilon\in\Qbb$ as follows:
$$P^\epsilon:=\{x\in\Qbb^n\,\mid\, Ax\geq B+\epsilon C\}.$$

We do not exclude the case where some lines of $A$ are zero. Note that $P^\epsilon$ can be empty (even for all $\epsilon\in\Qbb$). 

If there exists $\epsilon\in\Qbb$ such that $P^\epsilon$ is a (not empty) polytope (ie is bounded), then there is no non-zero $x\in\Qbb^n$ satisfying $Ax\geq 0$. Inversely, if there exists $\epsilon\in\Qbb$ such that $P^\epsilon$ is not bounded then there exists a non-zero $x\in\Qbb^n$ satisfying $Ax\geq 0$ (because $P^\epsilon$ contains at least an affine half-line and $x$ can be taken to be a generator of the direction of this half-line).

From now on, we suppose that there is no non-zero $x\in\Qbb^n$ satisfying $Ax\geq 0$.\\

Let $I_0:=\{1,\dots,m\}$.
For any matrix $\mathcal{M}$ and any $i\in I_0$, we denote by $\mathcal{M}_i$ the matrix consisting of the line $i$ of $\mathcal{M}$. And more generally, for any subset $I$ of $I_0$ we denote by $\mathcal{M}_I$ the matrix consisting of the lines $i\in I$ of $\mathcal{M}$.

Let $\epsilon\in\Qbb$. We denote by $\mathcal{H}_i^\epsilon$ 

the hyperplane $\{x\in\Qbb^n\,\mid\, A_ix=B+\epsilon C\}$. 
For any $I\subset I_0$, denote by $F_I^\epsilon$ the face of $P^\epsilon$ defined by $$F_I^\epsilon:=(\bigcap_{i\in I}\mathcal{H}_i^\epsilon)\cap Q^\epsilon.$$ Note that for any face $F^\epsilon$ of $P^\epsilon$ there exists a unique maximal $I\subset I_0$ such that $F^\epsilon=F_I^\epsilon$ (we include the empty face and $P^\epsilon$ itself). \\

Let $I\subset I_0$. Define $\Omega_{I,I_0}^0$ to be the set of $\epsilon\in\Qbb$ such that $F_I^\epsilon$ is not empty; define $\Omega_{I,I_0}^1$ to be the set of $\epsilon\in\Qbb$ such that, if $I'\subset I_0$ satisfies $F_I^\epsilon= F_{ I'}^\epsilon$, then $I'\subset I$.
In other words, $\epsilon\in\Omega_{I,I_0}^0$ if and only if there exists $x\in\Qbb^n$ such that $Ax\geq B+\epsilon C$, and $\epsilon\in\Omega_{I,I_0}^1$ if and only if there exists $x\in\Qbb^n$ such that $A_I=B_I+\epsilon C_I$ and $A_{I_0\backslash I}>B_{I_0\backslash I}+\epsilon C_{I_0\backslash I}$.

To make the notations not to heavy, we often write $i$ instead of $\{i\}$, for any $i\in I_0$.
Remark that if $\epsilon\in\Omega_{\emptyset,I_0}^1$, the polytope $P^\epsilon$ is of dimension~$n$ (ie has a non-empty interior). And, for any $i\in I_0$, if $\epsilon\in\Omega_{i,I_0}^1$ and $A_i\neq 0$, $F_i^\epsilon$ is a facet of $P^\epsilon$.\\

Now, we define an equivalence on subfamily of $(P^\epsilon)_{\epsilon\in\Qbb}$, that we extend to another family $(Q^\epsilon)_{\epsilon\in\Qbb}$ of polytopes constructed later (see Definition~\ref{defi:polytopeQ}), and corresponding to the equivalence of $G/H$-polytopes given in Definition~\ref{def:G/H-eq} restricted to the family used in Section~\ref{sec:MMP} (see Proposition~\ref{prop:eqeq}).

Let $K_0$ be the subset of $I_0$ consisting of indices of zero lines of $A$.

\begin{defi}\label{def:K-eq}
Let $K_0\subset K\subset I_0$.

We first define $\Omega_{K,I_0}^{max}:=\Omega_{\emptyset,I_0}^1\cap\bigcap_{i\in I_0\backslash K}\Omega^1_{i,I_0}$

Let $\epsilon$ and $\eta$ both in $\Omega_{K,I_0}^{max}$. 
We say that the polytopes $P^\epsilon$ and $P^\eta$ are equivalent if,
 for any $I\subset I_0$, $\epsilon$ and $\eta$ are either both in $\Omega_{I,I_0}^1$, either both in $\Qbb\backslash\Omega_{I,I_0}^0$, or both in $\Omega_{I,I_0}^0\backslash \Omega_{I,I_0}^1$. (In other words, $F_I^\epsilon$ and $F_I^\epsilon$ are either both not empty faces with $I$ maximal, either both empty, or both not empty with $I$ not maximal.)
\end{defi}

Some properties of the family $(P^\epsilon)_{\epsilon\in\Qbb}$ are listed in the following result.

\begin{teo}\label{th:mainpoly}
With the notations above, let $K_0\subset K\subset I_0$.
 \begin{enumerate}
 \item The set $\Omega_{K,I_0}^{max}$ is an open segment of $\Qbb$ (with extremities in $\Qbb\cup\{\pm\infty\}$, or empty).
 \item If $\Omega_{K,I_0}^{max}$ is not empty, there exist an integer $k$ and $\alpha_0<\cdots<\alpha_k$ in $\Qbb\cup\{\pm\infty\}$, such that the equivalent classes of polytopes in the family $(P^\epsilon)_{\epsilon\in\Omega_{K,I_0}^{max}}$ correspond exactly to the open segments $]\alpha_i,\alpha_{i+1}[$ for any $i\in\{0,\dots,k-1\}$ and the singletons $\{\alpha_i\}$ for any $i\in\{0,\dots,k\}$ with $\alpha_i\in\Qbb$.
 In particular, $\Omega_{K,I_0}^{max}=]\alpha_0,\alpha_k[$.
 \item Suppose that $\alpha_k\in\Qbb$. 
If $\alpha_k\in\Omega_{\emptyset,I_0}$, there exists $I_1\subset I_0$ containing $K$ such that $P^{\alpha_k}=P^{\alpha_k}_{I_1}:=\{x\in\Qbb^n\,\mid\,A_{I_1}x\geq B_{I_1}+\alpha_k C_{I_1}\}$ and $\alpha_k\in\Omega^{max}_{K,I_1}$ (defined for the family $(P^{\epsilon}_{I_1})_{\epsilon\in\Qbb}$).

If $\alpha_k\not\in\Omega_{\emptyset,I_0}$, there exist subsets $J_1$ and  $I_1$ of $I_0$ such that $J_1\cap I_1=\emptyset$, $P^{\alpha_k}=P^{\alpha_k}_{I_1,J_1}:=\{x\in\bigcap_{j\in J_1}\mathcal{H}_j^{\alpha_k}\,\mid\,A_{I_1}x\geq B_{I_1}+\alpha_k C_{I_1}\}$ and $\alpha_k$ is in $\Omega^{max}_{I_1\cap K,I_1}$ (defined for the family $(P^{\epsilon}_{I_1,J_1})_{\epsilon\in\Qbb}$ of polytopes in $\bigcap_{j\in J_1}\mathcal{H}_j^{\alpha_k}$). 

We have a similar result on the other extremity of $\Omega_{K,I_0}^{max}$ if $\alpha_0\in\Qbb$.
 \end{enumerate}
\end{teo} 

\subsection{Proof of Theorem~\ref{th:mainpoly}}

We study all the sets $\Omega_{I,I_0}^0$ and $\Omega_{I,I_0}^1$.

For any $I\subset I_0$, denote by $\bar{I}$ the complementary of $I$ in $I_0$.

\begin{lem}\label{lem:convexity}
For any $I\subset I_0$, the sets $\Omega_{I,I_0}^0$ and $\Omega_{I,I_0}^1$ are convex subsets of $\Qbb$ (may be empty).
\end{lem}

\begin{proof}
Let $\epsilon_1$ and $\epsilon_2$ be in $\Omega_{I,I_0}^0$. Then there exist $x_1$ and $x_2$ in $\Qbb^n$ such that $A_Ix_1=B_I+\epsilon_1C_I$, $A_Ix_2=B_I+\epsilon_2C_I$, $A_{\bar{I}}x_1\geq B_{\bar{I}}+\epsilon_1C_{\bar{I}}$ and $A_{\bar{I}}x_2\geq B_{\bar{I}}+\epsilon_2C_{\bar{I}}$. Then, for any rational number $\lambda\in [0,1]$, we get easily by adding equalities and inequalities that $A_I(\lambda x_1+(1-\lambda)x_2)=B_I+(\lambda\epsilon_1+(1-\lambda)\epsilon_2)C_I$ and $A_{\bar{I}}(\lambda x_1+(1-\lambda)x_2)\geq B_{\bar{I}}+(\lambda\epsilon_1+(1-\lambda)\epsilon_2)C_{\bar{I}}$, so that $\lambda\epsilon_1+(1-\lambda)\epsilon_2$ is in $\Omega_{I,I_0}^0$.

By replacing $\geq$ by $>$, we similarly prove the convexity of $\Omega_{I,I_0}^1$.
\end{proof}

\begin{lem}\label{lem:open}
For any $I\subset I_0$, the set $\Omega_{I,I_0}^1$ is either reduced to a point, or an open subset of $\Qbb$ (may be empty). Moreover, if $\Omega_{I,I_0}^1$ is reduced to a point $\epsilon_0$, the system $A_IX=B_I+\epsilon C_I$ has at least a solution only when $\epsilon=\epsilon_0$; in particular, in that case, $\Omega_{I,I_0}^0$ is also reduced to $\epsilon_0$.
\end{lem}

\begin{proof}
Suppose that $\Omega_{I,I_0}^1$ is neither empty nor reduced to a point. Let $\epsilon_1$ and $\epsilon_2$ two different points of $\Omega_{I,I_0}^1$. Then there exist $x_1$ and $x_2$ in $\Qbb^n$ such that $A_Ix_1=B_I+\epsilon_1C_I$, $A_Ix_2=B_I+\epsilon_2C_I$, $A_{\bar{I}}x_1> B_{\bar{I}}+\epsilon_1C_{\bar{I}}$ and $A_{\bar{I}}x_2> B_{\bar{I}}+\epsilon_2C_{\bar{I}}$. In particular both linear systems $A_IX=B_I$ and $A_IX=C_I$ have solutions, hence for any $\epsilon\in\Qbb$, the linear system $A_IX=B_I+\epsilon C_I$ has a set of solutions, depending continuously on $\epsilon$. More precisely, if $x_B$ and $x_C$ are some solutions of $A_IX=B_I$ and $A_IX=C_I$ respectively, then the set of solutions of  $A_IX=B_I+\epsilon C_I$ is $\ker(A_I)+x_B+\epsilon x_C$. Now, let $\epsilon\in\Omega_{I,I_0}^1$, and let $y:=x_0+x_B+\epsilon x_C$ with $x_0\in\ker(A_I)$ such that $A_{\bar{I}}y>B_{\bar{I}}+\epsilon C_{\bar{I}}$. Then there clearly exists a neighborhood $\mathcal{V}$ of $\epsilon$ in $\Qbb$ such that $A_{\bar{I}}(y+(\eta-\epsilon)x_C)>B_{\bar{I}}+\eta C_{\bar{I}}$ for all $\eta\in\mathcal{V}$. We also check easily that $A_I(y+(\eta-\epsilon)x_C)=B_I+\eta C_I$ so that $\mathcal{V}\subset\Omega_{I,I_0}^1$.

To prove the last statement, we suppose that the system $A_IX=B_I+\epsilon C_I$ has a solution for at least two values of $\epsilon$, and then by the same arguments as above, we deduce that $\Omega_{I,I_0}^1$ is open, which gives a contradiction. The result follows immediately.
\end{proof}

\begin{lem}\label{lem:gendarme}
For any $I\subset I_0$ such that $\Omega_{I,I_0}^1$ is not empty, we have  $\Omega_{I,I_0}^1\subset\Omega_{I,I_0}^0\subset\overline{\Omega_{I,I_0}^1}$.
\end{lem}
\begin{proof}
The first inclusion is obvious. 

To prove the second one, remark that for all $\epsilon_1\in\Omega_{I,I_0}^1$ and $\epsilon_2\in\Omega_{I,I_0}^0$, the segment $[\epsilon_1,\epsilon_2[$ (or $]\epsilon_2,\epsilon_1]$)  is contained in $\Omega_{I,I_0}^1$. Indeed, if there exists $x_1$ and $x_2$ such that $A_Ix_1=B_I+\epsilon_1C_I$, $A_Ix_2=B_I+\epsilon_2C_I$, $A_{\bar{I}}x_1> B_{\bar{I}}+\epsilon_1C_{\bar{I}}$ and $A_{\bar{I}}x_2\geq B_{\bar{I}}+\epsilon_2C_{\bar{I}}$, then for any $\lambda\in]0,1]$, $x_\lambda:=\lambda x_1+(1-\lambda)x_2$ and $\epsilon_\lambda:=\lambda \epsilon_1+(1-\lambda)\epsilon_2$ clearly satisfy $A_Ix_\lambda=B_I+\epsilon_\lambda C_I$, $A_{\bar{I}}x_\lambda> B_{\bar{I}}+\epsilon_\lambda C_{\bar{I}}$, so that $\epsilon_\lambda\in\Omega_{I,I_0}^1$. 
This remark implies directly that any $\epsilon_2\in\Omega_{I,I_0}^0$ is in $\overline{\Omega_{I,I_0}^1}$.
\end{proof}

\begin{lem}\label{lem:closed}
Let $I\subset I_0$ be such that $\Omega_{I,I_0}^1$ is not empty. Then the supremum and the infimum (well-defined in $\Rbb$) of $\Omega_{I,I_0}^1$ are either infinite or rational numbers. Moreover, $\Omega_{I,I_0}^0=\overline{\Omega_{I,I_0}^1}$.
\end{lem}

\begin{proof}
If $\Omega_{I,I_0}^1$ is reduced to a point, by Lemma~\ref{lem:gendarme} we have nothing to prove, so we suppose that $\Omega_{I,I_0}^1$ contains at least two points. 

Here (and only here), we need to consider the family of polytopes $(P^\epsilon)_{\epsilon\in\Rbb}$. Note that all definitions and all results given until now are still available replacing $\Qbb$ by $\Rbb$.
When it is necessary, we denote by $\Omega(\Qbb)$ and by $\Omega(\Rbb)$ the sets $\Omega$ defined at the beginning of the section respectively in $\Qbb$ and in $\Rbb$.
Now, by  Lemmas~\ref{lem:convexity} and~\ref{lem:gendarme}, it is enough to prove that, if they are finite, the supremum and the infimum of $\Omega_{I,I_0}^1(\Rbb)$ are rational numbers contained in the set $\Omega_{I,I_0}^0(\Qbb)$.

Suppose that the supremum $\epsilon_0\in\Rbb$ of $\Omega_{I,I_0}^1(\Rbb)$ is finite. Let $\epsilon_1\in\Omega_{I,I_0}^1(\Rbb)$. Then, for any $\epsilon\in[\epsilon_1,\epsilon_0]$, the polytope $P^\epsilon$ is contained in the polytope $$P^{[\epsilon_1,\epsilon_0]}:=\{x\in\Rbb^n\,\mid\,Ax\geq \Min(B+\epsilon_1C,B+\epsilon_0C)\},$$ where the mimimum $\Min$ is taken line by line.  
Now for all $\epsilon\in[\epsilon_1,\epsilon_0[$, let $x^\epsilon\in F_I^\epsilon$ such that  $A_{\bar{I}}x^\epsilon> B_{\bar{I}}+\epsilon C_{\bar{I}}$ (and $A_Ix^\epsilon=B_I+\epsilon C_I$). Since the points $x^\epsilon$ are in the compact set $P^{[\epsilon_1,\epsilon_0]}$ of $\Rbb^n$, there exists $x_0\in\Rbb^n$ such that $A_Ix_0=B_I+\epsilon_0C_I$ and $A_{\bar{I}}x_0\geq B_{\bar{I}}+\epsilon_0C_{\bar{I}}$, ie $x_0\in F_I^{\epsilon_0}$. It means that $\epsilon_0$ is in $\Omega_{I,I_0}^0(\Rbb)$.

Define the maximal subset $J$ of $I_0$ containing $I$ such that for any $x\in F_I^{\epsilon_0}$ we have $A_Jx=B_J+\epsilon_0C_J$. In particular, $A_{\bar{J}}x_0> B_{\bar{J}}+\epsilon_0C_{\bar{J}}$. We now prove that, for any $\epsilon_1\in\Omega_{I,I_0}^1$, the face $F_J^{\epsilon_1}$ is empty. Indeed, if it is not empty, there exists $x_1\in\Rbb^n$ such that $A_Jx_1=B_J+\epsilon_1C_J$ and $A_{\bar{J}}x_1\geq B_{\bar{J}}+\epsilon_1C_{\bar{J}}$. Let $\eta>0$, let $x_2:=x_0+\eta(x_0-x_1)$ and let $\epsilon_2:=\epsilon_0+\eta(\epsilon_0-\epsilon_1)$. Then we have $A_Jx_2=B_J+\epsilon_2C_J$ and, for $\eta$ small enough, we have $A_{\bar{J}}x_2> B_{\bar{J}}+\epsilon_2C_{\bar{J}}$.
Hence $F_I^{\epsilon_2}\supset F_J^{\epsilon_2}$ is not empty, which contradicts the fact that $\epsilon_2>\epsilon_0$ is not in $\Omega_{I,I_0}^0$ (see Lemma~\ref{lem:gendarme}). 

We now claim that the intersection of the vector space $\Im(A_J)$ with the affine line $\{B_J+\epsilon C_J\,\mid\,\epsilon\in\Rbb\}$ is reduced to the point $\epsilon_0$. Indeed, $\epsilon_0$ clearly belongs to this intersection, which can be either reduced to a point or the affine line. But, if some $\epsilon\in\Omega_{I,I_0}^1$ is in this intersection, then there exists $x_3\in\Rbb^n$ such that $A_Jx_3=B_J+\epsilon C_J$. Then, for $\eta>0$ small enough, we can prove that $\eta x_3+(1-\eta)x_0$ is in $F_J^{\epsilon_1}$ with $\epsilon_1=\eta\epsilon+(1-\eta)\epsilon_0$, which is not possible because $F_J^{\epsilon_1}$ is necessarily empty. To conclude that $\epsilon_0\in\Qbb$, it is now enough to recall that $A$, $B$ and $C$ have rational coefficients.

It remains to prove that $x_0$ can be chosen in $\Qbb^n$ so that $\epsilon_0$ is in $\Omega_{I,I_0}^0(\Qbb)$. But $x_0$ can be taken in an open set of the set of solutions of $A_Jx=B_J+\epsilon_0C_J$. We conclude by noticing that $B_J+\epsilon_0C_J$ have rational coefficients.

For the infimum of $\Omega_{I,I_0}^1$, it is the same proof (for example by replacing $\epsilon$ and $C$ by their opposites).
\end{proof}

\begin{cor}\label{cor:closed}
Let $I\subset I_0$. We have the following cases:
\begin{enumerate}
\item $\Omega_{I,I_0}^1$ is an open segment of $\Qbb$ (with extremities in $\Qbb\cup\{\pm\infty\}$) and $\Omega_{I,I_0}^0$ is the closed segment $\overline{\Omega_{I,I_0}^1}$ (may both be empty);
\item $\Omega_{I,I_0}^1$  is a reduced to a point $\epsilon_0$, then it equals $\Omega_{I,I_0}^0$ and the system $A_IX=B_I+\epsilon C_I$ has at least a solution only when $\epsilon=\epsilon_0$;
\item $\Omega_{I,I_0}^1$ is empty and $\Omega_{I,I_0}^0$ is the closure of some $\Omega_{J,I_0}^1$ with $I\varsubsetneq J\subset I_0$.
\end{enumerate}

Moreover, if $I=\emptyset$ or reduced to an index $i\in I_0\backslash K_0$, then only Cases~1 and 3 are possible, so that $\Omega_{I,I_0}^1$ is open (or empty).
\end{cor}

\begin{proof}
The two first items, and the fact that we cannot have more cases, can be deduced directly from Lemmas~\ref{lem:open} and \ref{lem:closed}.

Suppose now that $\Omega_{I,I_0}^1$ is empty, but not $\Omega_{I,I_0}^0$. 
Let $J\subset I_0$ be the maximal set containing $I$ such that, for all $\epsilon\in\Qbb$, $F_I^\epsilon= F_J^\epsilon$. Then $J\neq I$ by hypothesis. Indeed, if $I=J$ then there exist subsets $J_1,\dots,J_k$ of $I_0$ whose intersection is $I$, and also rational numbers $\epsilon_1,\dots,\epsilon_k$ respectively in $\Omega_{J_1,I_0}^1,\dots,\Omega_{J_k,I_0}^1$. Hence $\frac{\epsilon_1+\cdots+\epsilon_k}{k}$ is in $\Omega_{I,I_0}^1$ that is not possible. 

Then $\Omega_{I,I_0}^0=\Omega_{J,I_0}^0$ and by maximality of $J$, $\Omega_{J,I_0}^1$ is not empty. We conclude by using the first two items.

Now, for the last statement remark that, if $I=\emptyset$ or reduced to an index $i\in I_0\backslash K_0$, the system $A_IX=B_I+\epsilon C_I$ has solutions for all $\epsilon$ and then, we cannot be in Case~2.
\end{proof}

We now prove three lemmas to get the connexity of $K$-equivalent classes of polytopes.

\begin{lem}\label{lem:separationextremities}
Let $I\subset I_0$ such that $\Omega_{I,I_0}^0$ is bounded, but neither empty nor reduced to a point. Denote by  $\epsilon_1$ and $\epsilon_2$ the extremities of $\Omega_{I,I_0}^0$. Then, there exist $J_1$ and $J_2$ containing $I$ such that $\{\epsilon_1\}=\Omega_{J_1,I_0}^1$ and $\{\epsilon_2\}=\Omega_{J_2,I_0}^1$. In particular, if $\epsilon_1$ and $\epsilon_2$ are in $\Omega^{max}_{K,I_0}$, $P^{\epsilon_1}$ is not equivalent to $P^{\epsilon_2}$.
\end{lem}

\begin{proof}
By Corollary \ref{cor:closed}, we can suppose that $\Omega_{I,I_0}^0=\overline{\Omega_{I,I_0}^1}=[\epsilon_1,\epsilon_2]$. Then, there exists $j\in\bar{I}$ such that $x\in Q^{\epsilon_1}$ and $A_Ix=B_I+\epsilon_1 C_I$ imply that $A_jx=B_j+\epsilon C_j$. Let $J_1$ be the union of $I$ with the set of all such indices $j$. In particular, $\epsilon_1\in\Omega_{J_1,I_0}^1$. 
Moreover, since $\Omega_{J_1,I_0}^1$ is contained in $\Omega_{I,I_0}^0=[\epsilon_1,\epsilon_2]$, it cannot be open, and then by Corollary~\ref{cor:closed}, $\Omega_{J_1,I_0}^1$ is reduced to $\epsilon_1$.

Similarly, we prove the existence of $J_2$.

In particular, if $\epsilon_1$ and $\epsilon_2$ are in $\Omega^{max}_{K,I_0}$, by definition, the classes of $P^{\epsilon_1}$ and $P^{\epsilon_2}$ are both reduced to one polytope and then distinct.
\end{proof}

\begin{lem}\label{lem:pointtosegement}
Let $I\subset I_0$ such that $\Omega_{I,I_0}^0=\{\epsilon_0\}$, with $\epsilon_0\in\Omega_{K,I_0}^{max}$. There exists $I_1\subset I_0$ such that the upper extremity of $\Omega_{I_1,I_0}^0$ is $\epsilon_0$. 
\end{lem}

\begin{proof}
First, we can suppose that $I$ is the maximal subset of $I_0$ such that $\Omega_{I,I_0}^0=\{\epsilon_0\}$, (ie such that $\Omega_{I,I_0}^1=\{\epsilon_0\}$). Then there exists $x_0\in\Qbb^n$ satisfying $A_Ix_0=B_I+\epsilon_0C_I$ and $A_{\bar{I}}x_0>B_{\bar{I}}+\epsilon_0C_{\bar{I}}$.

Consider a subset $I_1$ of $I$ such that there exists $\epsilon_1<\epsilon_0$ in $\Omega_{I_1,I_0}^0$. Choose $I_1$ maximal in $I$ with this property. It exists (but can be empty) because $\epsilon_0$ is in the open set $\Omega_{K,I_0}^{max}$. Moreover $I_1\neq I$.
Then there exist $\epsilon_1<\epsilon_0$ and $x_1\in \Qbb^n$ such that $A_{I_1}x_1=B_{I_1}+\epsilon_1 C_{I_1}$ and $A_{\bar{I_1}}x_1\geq B_{\bar{I_1}}+\epsilon_1 C_{\bar{I_1}}$. Hence, for any $\epsilon\in]\epsilon_1,\epsilon_0[$, the element $x^\epsilon:=\frac{(\epsilon_0-\epsilon)x_1+(\epsilon-\epsilon_1)x_0}{\epsilon_0-\epsilon_1}$ satisfies $A_{I_1}x^\epsilon=B_{I_1}+\epsilon C_{I_1}$ and $A_{\bar{I_1}}x^\epsilon> B_{\bar{I_1}}+\epsilon C_{\bar{I_1}}$ by maximality of $I_1$. It proves that  $]\epsilon_1,\epsilon_0[\subset\Omega_{I_1,I_0}^1$. In particular, with Corollary~\ref{cor:closed}, $\Omega_{I_1,I_0}^1$ is open with closure equal to  $\Omega_{I_1,I_0}^0$. 

It is now enough to prove that the upper extremity of $\Omega_{I_1,I_0}^0$  is $\epsilon_0$. Suppose the contrary, then there exists $\epsilon_2>\epsilon_0$ and $x_2\in\Qbb^n$ such that $A_{I_1}x_2=B_{I_1}+\epsilon_2 C_{I_1}$ and $A_{\bar{I_1}}x_2\geq B_{\bar{I_1}}+\epsilon_2 C_{\bar{I_1}}$.
For any $\epsilon\in]\epsilon_0,\epsilon_2[$, the element $y^\epsilon:=\frac{(\epsilon_2-\epsilon)x_0+(\epsilon-\epsilon_0)x_2}{\epsilon_2-\epsilon_0}$ satisfies $A_{I_1}y^\epsilon=B_{I_1}+\epsilon C_{I_1}$ and $A_{\bar{I_1}}y^\epsilon> B_{\bar{I_1}}+\epsilon C_{\bar{I_1}}$. Moreover, for $\epsilon=\epsilon_0$ (ie $y^\epsilon=x_0$), we have $A_{I}y^\epsilon=B_{I}+\epsilon C_{I}$ and $A_{\bar{I}}y^\epsilon> B_{\bar{I}}+\epsilon C_{\bar{I}}$. Then by continuity, for any $\epsilon<\epsilon_0$ big enough, we have $A_{I_1}y^\epsilon=B_{I_1}+\epsilon C_{I_1}$, $A_{\bar{I}}y^\epsilon> B_{\bar{I}}+\epsilon C_{\bar{I}}$ and $A_{I\backslash I_1}y^\epsilon<B_{I\backslash I_1}+\epsilon C_{I\backslash I_1}$. Choose such an $\epsilon$ in $]\epsilon_1,\epsilon_0[$. Hence, there exist $z^\epsilon$ in the segment 
$]x^\epsilon,y^\epsilon[$ and some $j\in I\backslash I_1$ such that $A_{I_1}z^\epsilon=B_{I_1}+\epsilon C_{I_1}$, $A_{\bar{I_1}}z^\epsilon\geq B_{\bar{I_1}}+\epsilon C_{\bar{I_1}}$ and $A_{j}z^\epsilon=B_{j}+\epsilon C_{j}$. This contradicts the maximality of $I_1$, so that we have proved that the upper extremity of $\Omega_{I_1,I_0}^0$ is $\epsilon_0$.
\end{proof}

\begin{lem}\label{lem:connexity}
Let $\eta<\epsilon$ in $\Omega^{max}_{K,I_0}$ such that there exists $I\subset I_0$ satisfying $\Omega_{I,I_0}^0\subset]\eta,\epsilon[$.
Then $P^{\eta}$ is not equivalent to $P^{\epsilon}$.
\end{lem}

\begin{proof}
Among the set of subsets $I$ of $I_0$ satisfying $\Omega_{I,I_0}^0\subset]\eta,\epsilon[$, choose the one whose lower extremity is minimal. 

There are two cases: in the first one, $\Omega_{I,I_0}^0=\{\epsilon_0\}$.
By the Lemma~\ref{lem:pointtosegement}, there exist $I_1\subset I_0$ and $\epsilon_1<\epsilon_0$ (may be $-\infty$) such that $\Omega_{I_1,I_0}^0=[\epsilon_1,\epsilon_0]$. By hypothesis of minimality, $\eta$ is in $\Omega_{I_1,I_0}^0$. But $\epsilon$ is not in $\Omega_{I_1,I_0}^0$, so that $P^{\eta}$ is not equivalent to $P^{\epsilon}$.

Consider the second case:  $\Omega_{I,I_0}^0=[\epsilon_1,\epsilon_2]$ with $\epsilon_1<\epsilon_2$.
Let $J$ be the maximal subset of $I_0$ such that  $x\in Q^{\epsilon_1}$ and $A_Ix=B_I+\epsilon_1 C_I$ imply that $A_Jx=B_J+\epsilon_1 C_J$. In particular, $J$ strictly contains $I$ and $\epsilon_1\in\Omega_{J,I_0}^1$.
Then $\Omega_{J,I_0}^1\subset \Omega_{I,I_0}^0$ cannot be open, hence it is reduced to a point and we conclude by the first case.
\end{proof}

\begin{proof}[Proof of Theorem~\ref{th:mainpoly}]
\begin{enumerate}
\item
It comes directly from the last statement of Corollary~\ref{cor:closed}.

\item Now, let $\epsilon_0\in\Omega_{K,I_0}^{max}$. Prove that the equivalent class of $P^{\epsilon_0}$ corresponds to a segment of $\Omega_{K,I_0}^{max}$ or is reduced to $P^{\epsilon_0}$.

If there exists $I\subset I_0$ such that $\Omega_{I,I_0}^1=\{\epsilon_0\}$, then the equivalent class of $P^{\epsilon_0}$ is clearly reduced to $P^{\epsilon_0}$.
By Lemma~\ref{lem:separationextremities}, we get the same conclusion if there exists $I\subset I_0$ such that $\Omega_{I,I_0}^1$ is not empty and $\epsilon_0\in\Omega_{I,I_0}^0\backslash\Omega_{I,I_0}^1$ (ie if there exists $I\subset I_0$ such that $\epsilon_0$ is an extremity of $\Omega_{I,I_0}^0$, by Corollary~\ref{cor:closed} Case~3).

Suppose now that for all subsets $I$ of $I_0$, $\Omega_{I,I_0}^1\neq\{\epsilon_0\}$, and $\epsilon_0$ is not an extremity of $\Omega_{I,I_0}^0$. Then the set of $\epsilon\in\Omega_{K,I_0}^{max}$ such that $P^\epsilon$ is equivalent to $P^{\epsilon_0}$ is the intersection of open segment of $\Qbb$. Indeed, it is the intersection of open segments of type $\Omega_{I,I_0}^1$ and connected components of some $\Qbb\backslash\Omega_{I,I_0}^0$ (by Lemma~\ref{lem:connexity}).

\item Now, let $\epsilon_1$ be an extremity of $\Omega_{K,I_0}^{max}$. Suppose that $\epsilon_1\in\Omega_{\emptyset,I_0}^1$. Define $I_2$ to be the set of indices $i\in I_0\backslash K$ such that $\epsilon_1$ is not in $\Omega_{i,I_0}^1$. Denote by $I_1$ the complementary of $I_2$ in $I_0$. Then by definition, $I_1$ contains $K$ and $\epsilon_1$ is in the set $\Omega_{K,I_1}^{max}$. Denote by $P_{I_1}^{\epsilon_1}$ the polyhedron $\{x\in\Qbb^n\,\mid\,A_{I_1}x\geq B_{I_1}+\epsilon_1 C_{I_1}\}$, and prove that $P^{\epsilon_1}=P_{I_1}^{\epsilon_1}$. We clearly have $P^{\epsilon_1}\subset P_{I_1}^{\epsilon_1}$, so that we have to prove that $P_{I_1}^{\epsilon_1}$ is contained in all closed half-spaces $\mathcal{H}_i^{+\epsilon_1}:=\{x\in\Qbb^n\,\mid\,A_ix\geq B_i+\epsilon_1 C_i\}$ with $i\in I_2$, or equivalently that the interior $\{x\in\Qbb^n\,\mid\,A_{I_1}x>B_{I_1}+\epsilon_1 C_{I_1}\}$ of $P_{I_1}^{\epsilon_1}$ (that is not empty because $\epsilon_1\in\Omega_{\emptyset,I_0}^1$) is contained in all open half-spaces $\mathcal{H}_i^{++\epsilon_1}:=\{x\in\Qbb^n\,\mid\,A_ix> B_i+\epsilon_1 C_i\}$ with $i\in I_2$. Since the interior of $P_{I_1}^{\epsilon_1}$ intersects the connected component $\bigcap_{i\in I_2}\mathcal{H}_i^{++\epsilon_1}$ of $\Qbb^n\backslash\bigcup_{i\in I_2}\mathcal{H}_i^{\epsilon_1}$, it is enough to prove that the interior of $P_{I_1}^{\epsilon_1}$ does not intersect the union of hyperplanes $\bigcup_{i\in I_2}\mathcal{H}_i^{\epsilon_1}$.

Let $i\in I_2$. Suppose that the interior of $P_{I_1}^{\epsilon_1}$ intersects the hyperplane $\mathcal{H}_i^{\epsilon_1}$. There exists $x_1\in\Qbb^n$ such that $A_{I_1}x_1> B_{I_1}+\epsilon_1 C_{I_1}$ and $A_ix_1=B_i+\epsilon_1C_i$. Since there also exists $x_0\in\Qbb^n$ such that $Ax_0>B+\epsilon_1C$, an element $x_2$ of the segment $[x_1,x_0[$ is in the interior of $P_{I_1}^{\epsilon_1}$, in $P^{\epsilon_1}$ and in $\mathcal{H}_j^{\epsilon_1}$ for at least one $j\in I_2$. Let $J\subset I_2$ the maximal set such that $A_Jx_2=B_J+\epsilon_1C_J$. Then the set $\mathcal{C}:=\{x\in\Qbb^n\,\mid\,A_Jx\geq B_J+\epsilon_1C_J\}$ is a cone (of apex $x_2+\ker(A_J)$) with non-empty interior. Let $j\in J$ such that $\mathcal{H}_j^{\epsilon_1}\cap \mathcal{C}$ is a facet of $\mathcal{C}$.  By Lemma~\ref{lem:nocolinearlines} below, all hyperplanes $\mathcal{H}_j^{\epsilon_1}$ with $j\in J$ are different, then there exists $y\in\mathcal{H}_j^{\epsilon_1}\cap\mathcal{C}$ that is not in the other hyperplanes $\mathcal{H}_i^{\epsilon_1}$ with $j\neq i\in J$. Since $x_2$ is in the interior of $P_{I_0\backslash J}^{\epsilon_1}$, for any rational number $\lambda>0$ small enough, the point $x_2+\lambda(y-x_2)$ is in the interior of $P_{I_0\backslash J}^{\epsilon_1}$, also in $\mathcal{H}_j^{\epsilon_1}\cap\mathcal{C}$ and not in $\mathcal{H}_i^{\epsilon_1}$ for all $i\in J$ different from $j$. This contradicts the fact that $\epsilon_1\not\in\Omega_{j,I_0}^1$.

Suppose now that $\epsilon_1\not\in\Omega_{\emptyset,I_0}^1$.
Let $J_1$ the maximal subset of $I_0$ such that $P^{\epsilon_1}$ is contained in the subspace $\mathcal{H}_{J_1}:=\bigcap_{j\in J_1}\mathcal{H}_j^{\epsilon_1}$ of $\Qbb^n$.
Then $P^{\epsilon_1}$ is clearly the polytope $\{x\in\bigcap_{j\in J_1}\mathcal{H}_j^{\epsilon_1}\,\mid\,A_{I_0\backslash J_1}x\geq B_{I_0\backslash J_1}+\epsilon_1 C_{I_0\backslash J_1}\}$.
Then $\epsilon_1\in\Omega_{\emptyset,I_0\backslash J_1}^1$ (defined in for the family of polytopes in  $\mathcal{H}_{J_1}$). If $\epsilon_1\in\Omega_{K\cap(I_0\backslash J_1),I_0\backslash J_1}^{max}$, then $I_1=I_0\backslash J_1$ gives the result. If not, $\epsilon_1$ is an extremity of $\Omega_{K\cap(I_0\backslash J_1),I_0\backslash J_1}^{max}$ and $\epsilon_1\in\Omega_{\emptyset,I_0\backslash J_1}^1$. Then we apply that we prove in the previous paragraph, in order to find $I_1\subset I_0\backslash J_1$ that gives the result.
\end{enumerate}
\end{proof}

\begin{lem}\label{lem:nocolinearlines}
If $\Omega_{K,I_0}^{max}$ is not empty, then for any $i$ and $j$ in $I_0\backslash K$, the existence of $\lambda\in\Qbb_{>0}$ such that $A_i=\lambda A_j$, implies that $i=j$.
\end{lem}

\begin{proof}
Let $\epsilon\in\Omega_{K,I_1}^{max}$. Suppose that there exists $i\neq j$ in $I_0\backslash K$ such that $A_i=\lambda A_j$ for some $\lambda\in\Qbb_{>0}$. Then, since $\epsilon\in\Omega_{i,I_0}^1\cap\Omega_{j,I_0}^1$, there exist $x$ and $y$ in $\Qbb^n$, such that $A_ix=B_i+\epsilon C_i$, $A_jx>B_j+\epsilon C_j$, $A_jy=B_j+\epsilon C_j$ and $A_iy>B_i+\epsilon C_i$.
And we have $\lambda(B_j+\epsilon C_j)<\lambda A_jx=A_ix=B_i\epsilon C_i<A_iy=\lambda A_jy=\lambda(B_j+\epsilon C_j)$, which gives a contradiction.
\end{proof}

\subsection{A second one-parameter family of polytopes}

In section~\ref{sec:MMP}, we apply Theorem~\ref{th:mainpoly} to a family $(\tilde{Q}^\epsilon)_{\epsilon\in\Qbb_{\geq 0}}$ constructed by iteration from the family $(P^\epsilon)_{\epsilon\in\Qbb}$, using the last statement of Theorem~\ref{th:mainpoly} in the following way.

\begin{defi}\label{defi:polytopeQ}
Let $A$, $B$ and $C$ be as in Theorem~\ref{th:mainpoly}. Let $K\subset I_0$ containing $K_0$.
Suppose that $0\in\Omega_{K,I_0}^{max}$.
For any $\epsilon\in\Qbb_{\geq 0}\cap\Omega_{K,I_0}^{max}$, define $\tilde{Q}^\epsilon$ to be $P^\epsilon$.
If $\Omega_{K,I_0}^{max}$ has a finite supremum $\epsilon_1$ and if $\epsilon_1\in\Omega_{\emptyset,I_0}$ consider $I_1$ as in Theorem~\ref{th:mainpoly}. Then define for any $\epsilon\in[\epsilon_1,+\infty[\cap\Omega^{max}_{K,I_1}$, $\tilde{Q}^\epsilon$ to be the polytope $P^\epsilon_{I_1}$.

If $\epsilon_1\not\in\Omega_{\emptyset,I_0}$, we stop the construction and we define $\tilde{Q}^\epsilon=\emptyset$ for any $\epsilon>\epsilon_1$ by convention.

Iterating the construction, we obtain a family of polytopes $(\tilde{Q}^\epsilon)_{\epsilon\in\Qbb_{\geq 0}}$.
Remark that the construction depends on $K$.

We say that $\tilde{Q}^\epsilon$ is $K$-equivalent to $\tilde{Q}^\eta$ in the family $(\tilde{Q}^\epsilon)_{\epsilon\in\Qbb_{\geq 0}}$ if they are both defined at the same step of the construction above, ie if they correspond to $P^\epsilon_I$ and $P^\eta_I$, with $\epsilon$ and $\eta$ both in $\Omega^{max}_{K,I}$ for some $I\subset I_0$ containing $K$, and $P^\epsilon_I$ and $P^\eta_I$ are equivalent (according to Definition~\ref{def:K-eq}).
\end{defi}

Applying Theorem~\ref{th:mainpoly} to each subfamilies $(P_{I_i}^\epsilon)$ of the family $(\tilde{Q}^\epsilon)_{\epsilon\in\Qbb_{\geq 0}}$, we obtain immediately the following result.

\begin{cor}[of Theorem~\ref{th:mainpoly}]\label{cor:mainpoly}
Suppose that $\tilde{Q}^0$ and $\tilde{Q}^\epsilon$ are equivalent for $\epsilon>0$ small enough.

There exist non-negative integers $k,j_0,\dots,j_k$, rational numbers $\alpha_{i,j}$ for $i\in\{0,\dots,k\}$ and $j\in\{0,\dots,j_i\}$ and $\alpha_{k,j_k+1}\in\Qbb_{>0}\cup\{+\infty\}$ ordered as follows with the convention that $\alpha_{i,j_i+1}=\alpha_{i+1,0}$ for any $i\in\{0,\dots,k-1\}$:
\begin{enumerate}
\item $\alpha_{0,0}=0$;
\item for any $i\in\{0,\dots,k\}$, and for any $j<j'$ in $\{0,\dots,j_i+1\}$ we have $\alpha_{i,j}<\alpha_{i,j'}$;
\end{enumerate}
and such that the $K$-equivalent classes in the family of polytopes $(\tilde{Q}^\epsilon)_{\epsilon\in\Qbb_{\geq 0}}$ is given by the following segments:

\begin{enumerate}
\item $[\alpha_{i,0},\alpha_{i,1}[$, with $i\in\{0,\dots,k\}$;
\item $]\alpha_{i,j},\alpha_{i,j+1}[$, with $i\in\{0,\dots,k\}$ and $j\in\{1,\dots,j_i\}$;
\item $\{\alpha_{i,j}\}$ with $i\in\{0,\dots,k\}$ and $j\in\{1,\dots,j_i\}$;
\item  if $\alpha_{k,j_k+1}\neq +\infty$, $\{\alpha_{k,j_k+1}\}$ and $]\alpha_{k,j_k},+\infty[$.
 \end{enumerate}
\end{cor}

\begin{rems}\label{rem:epsilonmax}
For any $i\in\{0,\dots,k\}$, the rational numbers $\alpha_{i,0},\dots,\alpha_{i,j_i}$ can come from several consecutive steps of the construction of $(\tilde{Q}^\epsilon)_{\epsilon\in\Qbb_{\geq 0}}$. Indeed, for example, $\alpha_{0,j_0}$ is not necessarily the supremum $\epsilon_1$ of $\Omega^{max}_{K,I_0}$ because $\{\tilde{Q}^{\epsilon_1}\}$ may be an equivalent class.

We have $\alpha_{k,j_k+1}= +\infty$ if, in the construction of $(\tilde{Q}^\epsilon)_{\epsilon\in\Qbb_{\geq 0}}$, some set $\Omega^{max}_{K,I}$ has no upper bound. And if $\alpha_{k,j_k+1}\neq +\infty$, it is the upper extremity of some $\Omega^{max}_{K,I}$ that is not an element of $\Omega_{\emptyset,I}^1$. In this latter case, $]\alpha_{k,j_k+1},+\infty[$ is the class of empty polytopes and $\alpha_{k,j_k+1}$ corresponds to the last not empty polytope. If $C>0$, we are always in that case.
\end{rems}

\section{MMP via a one-parameter family of polytopes} \label{sec:MMP}

Let $X$ be a projective horospherical $G$-variety, with open $G$-orbit isomorphic to $G/H$. Let $D$ be an ample $\Qbb$-Cartier divisor. Suppose that $X$ is $\Qbb$-Gorenstein (ie the canonical divisor $K_X$ of $X$ is $\Qbb$-Cartier). We keep the notations given in Sections~\ref{sec:horopoly}.

\subsection{The one-parameter family of polytopes}\label{sec:TheFamily}

To construct the one-parameter family of polytopes that permits to run the MMP from $X$, we start the same way as in the classical approach of the MMP. Indeed, for $\epsilon>0$ small enough, the divisor $D+\epsilon K_X$ is still ample (and $\Qbb$-Cartier by hypothesis), so that it defines a moment polytope $Q^\epsilon$ and a pseudo-moment polytope $\tilde{Q}^\epsilon$. More precisely, (for $\epsilon$ small enough) $Q^\epsilon:=\{x\in v^\epsilon+M_\Qbb\,\mid\, Ax\geq B+\epsilon C\}$ and $\tilde{Q}^\epsilon:=\{x\in M_\Qbb\,\mid\, Ax\geq \tilde{B}+\epsilon \tilde{C}\}$ where 
the matrices $A$, $B$, $C$, $\tilde{B}$, $\tilde{C}$ and the vector $v^\epsilon$ are defined below. 

Recall that $x_1,\dots,x_m$ denote the primitive elements of the rays of the colored fan of $X$ that are not generated by a vector $\alpha^\vee_M$ with $\alpha\in\mathcal{F}_X$. We choose an order in $S\backslash R$ and we then denote by $\alpha_1,\dots,\alpha_r$ its elements. We fix a basis $\mathcal{B}$ of $M$. Now define $A\in\mathcal{M}_{m+r,n}(\Qbb)$ whose first $m$ lines are the coordinates of the vectors $x_i$ in the basis $\mathcal{B}$ with $i\in\{1,\dots,m\}$ and whose last $r$ lines are the coordinates of the vectors $\alpha^\vee_{jM}$ in $\mathcal{B}$ with $j\in\{1,\dots,r\}$. Let $\tilde{B}$ be the column vector such that the pseudo-moment polytope of $D$ is defined by $\{x\in M_\Qbb\,\mid\, Ax\geq \tilde{B}\}$. In fact, if $D=\sum_{i=1}^m b_iX_i+\sum_{\alpha\in S\backslash R}b_\alpha D_\alpha$, then $\tilde{B}$ is the column matrix associated to the vector $(-b_1,\dots,-b_m,-b_{\alpha_1},\dots,-b_{\alpha_r})$. Similarly, the column matrix $\tilde{C}$ corresponds to the vector $(1,\dots,1,c_{\alpha_1},\dots,c_{\alpha_r})$, where $-K_X=\sum_{i=1}^mX_i+\sum_{\alpha\in S\backslash R}c_\alpha D_\alpha$ (the coefficients are explicitly defined with $c_\alpha=\langle 2\rho_P,\alpha^\vee\rangle$, where $\rho_P$ is the sum of positive roots of $G$ that are not roots of $P$). Now, define $v^\epsilon:=\sum_{\alpha\in S\backslash R}(b_\alpha-\epsilon c_\alpha)\varpi_\alpha$. Since $Q^\epsilon=v^\epsilon+\tilde{Q}^\epsilon$, we compute that $B$ and $C$ are respectively the column matrices associated to $(-b_1+\langle x_1,\sum_{\alpha\in S\backslash R}b_\alpha\varpi_\alpha\rangle,\dots,-b_m+\langle x_m,\sum_{\alpha\in S\backslash R}b_\alpha\varpi_\alpha\rangle,0,\dots,0)$ and $(1-\langle x_1,\sum_{\alpha\in S\backslash R}c_\alpha\varpi_\alpha\rangle,\dots,1-\langle x_m,\sum_{\alpha\in S\backslash R}c_\alpha\varpi_\alpha\rangle,0,\dots,0)$.

From the matrices $A$, $\tilde{B}$ and $\tilde{C}$ we construct the family $(\tilde{Q}^\epsilon)_{\epsilon\in\Qbb_{\geq 0}}$ of polytopes in $M_\Qbb$ as in Section~\ref{sec:poly}. And for any $\epsilon\geq 0$ we define $Q^\epsilon=v^\epsilon+\tilde{Q}^\epsilon$. Since $\tilde{C}>0$, by Remark~\ref{rem:epsilonmax}, there exists $\epsilon_{max}\in\Qbb_{>0}$ (it is the number $\alpha_{k,j_k+1}$ given in Corollary~\ref{cor:mainpoly}) such that for all $\epsilon\in[0,\epsilon_{max}[$, $Q^\epsilon$ is a $G/H$-polytope and for all $\epsilon>\epsilon_{max}$, $Q^\epsilon=\emptyset$. The polytope $Q^{\epsilon_{max}}$ is neither empty nor a $G/H$-polytope, but it is a $G/H^1$-poltyope with some subgroup $H^1$ of $G$ containing $H$ (see later).

\begin{prop}\label{prop:eqeq}
The two partitions of $[0,\epsilon_{max}[$ given by equivalence classes of $G/H$-polytopes (Definition~\ref{def:G/H-eq}) in the family $(Q^\epsilon)_{\epsilon\in[0,\epsilon_{max}[}$ and by equivalent classes in the family $(\tilde{Q}^\epsilon)_{\epsilon\in[0,\epsilon_{max}[}$  according to Definition~\ref{defi:polytopeQ} with $K=\{m+1,\dots,m+r\}$ are the same.
\end{prop}

\begin{proof}
Let $\epsilon$ and $\eta$ be two rational numbers. 

Suppose that $Q^\epsilon$ and $Q^\eta$ are equivalent $G/H$-polytopes. In Definition~\ref{def:G/H-eq}, take $j$ to be minimal, so that the hyperplanes $\mathcal{H}_1,\dots,\mathcal{H}_j$ and $\mathcal{H'}_1,\dots,\mathcal{H'}_j$ corresponds respectively to $\mathcal{H}_{i_1}^\epsilon,\dots,\mathcal{H}_{i_j}^\epsilon$ and $\mathcal{H}_{i_1}^\eta,\dots,\mathcal{H}_{i_j}^\eta$ for $J:=\{i_1,\dots,i_j\}\subset\{1,\dots,m+r\}$.
Then, if we denote by $I$ the union $J\cup K$, the polytopes $\tilde{Q}^\epsilon$ and $\tilde{Q}^\eta$ respectively equal $P_I^\epsilon$ and $P_I^\eta$ such that $\epsilon$ and $\eta$ are in $\Omega_{K,I}^{max}$. Then we have to prove that $P_{I}^\epsilon$ and $P_{I}^\eta$ are equivalent according to Definition~\ref{def:K-eq}, which comes directly from the definition of equivalence of $G/H$-polytopes.

Suppose now that $\tilde{Q}^\epsilon$ and $\tilde{Q}^\eta$ are equivalent according to Definition~\ref{defi:polytopeQ}. They are constructed at the same step so that there exists $I\subset\{1,\dots,m+r\}$ containing $K$ such that $\tilde{Q}^\epsilon$ and $\tilde{Q}^\eta$ respectively equal $P_I^\epsilon$ and $P_I^\eta$, with $\epsilon$ and $\eta$ in $\Omega_{K,I}^{max}$. Then $\tilde{Q}^\epsilon=\{x\in M_\Qbb\,\mid\, A_I\geq \tilde{B}_I+\epsilon \tilde{C}_I\}$ and $\tilde{Q}^\eta=\{x\in M_\Qbb\,\mid\, A_I\geq \tilde{B}_I+\eta \tilde{C}_I\}$).  This directly gives the first item of Definition~\ref{def:G/H-eq}. 
The second item is also clear from Definition~\ref{def:K-eq}. And the third item comes from the fact that $\Omega_{i,I}^0$ contains both $\epsilon$ and $\eta$ or no of the two, for any $i\in K$.
\end{proof}

\subsection{Construction of varieties and morphisms}\label{sec:constructionvarietes}

We apply Corollary~\ref{cor:mainpoly} to the family $(\tilde{Q}^\epsilon)_{\epsilon\in\Qbb_{\geq 0}}$. 
Then the family $(Q^\epsilon)_{\epsilon\in\Qbb_{\geq 0}}$ gives a list of $G/H$-embeddings:
\begin{enumerate}
\item $X_{i,j}$ for any  $i\in\{0,\dots,k\}$ and $j\in\{0,\dots,j_i\}$, respectively associated to moment polytopes $Q^\epsilon$ with $\epsilon\in]\alpha_{i,j},\alpha_{i,j+1}[$;
\item $Y_{i,j}$ for any $i\in\{0,\dots,k\}$ and $j\in\{1,\dots,j_i\}$, respectively associated to moment polytopes $Q^{\alpha_{i,j}}$;
\end{enumerate}

It also gives a projective horospherical $G$-variety $Z$ associated to the moment polytope $Q^{\alpha_{k,j_k+1}}=Q^{\alpha_{max}}$. Indeed,  let $M^1_\Qbb$ be the minimal vector subspace containing $\tilde{Q}^{\alpha_{max}}$ and let $M^1:=M^1_\Qbb\cap M$. Let $R^1$ be the union of $R$ with the set of $\alpha\in S\backslash R$ such that $Q^{\alpha_{max}}$ is contained in the wall $W_{\alpha,P}$. Then we define the subgroup $H^1$ of $P^1:=P_{R^1}$ to be  the intersection of kernels of characters of $P^1$ in $M^1$. Then $Q^{\alpha_{max}}$ is a $G/H^1$-polytope and corresponds to a $G/H^1$-embedding $Z$. Remark that, by definition $M^1\subset M$ and $R\subset R^1$ so that we have a projection $\pi:G/H\longrightarrow G/H^1$.\\

Now, by Proposition~\ref{prop:morph}, we get $G$-equivariant morphisms: 
\begin{enumerate} \item $\phi_{i,j}:X_{i,j-1}\longrightarrow Y_{i,j}$ for any $i\in\{0,\dots,k\}$ and $j\in\{1,\dots,j_i\}$;
\item $\phi_{i,j}^+:X_{i,j}\longrightarrow Y_{i,j}$ for any $i\in\{0,\dots,k\}$ and $j\in\{1,\dots,j_i\}$; \item $\phi_i:X_{i,j_i}\longrightarrow X_{i+1,0}$ for any $i\in\{0,\dots,k-1\}$;
\item and $\phi:X_{k,j_k}\longrightarrow Z$.
\end{enumerate}
In the next section, we prove that the morphisms $\phi_{i,j}$ and $\phi_{i,j}^+$ give flips (may be divisorial, see Remark~\ref{rem:epsilonmax} and Example~\ref{ex:horo5}), the morphisms $\phi_i$ are divisorial contractions and the morphism $\phi$ is a Mori fibration.

\subsection{Description of the contracted curves}

\begin{prop}
\begin{enumerate} \item For any $i\in\{0,\dots,k\}$ and $j\in\{1,\dots,j_i\}$, the curves $C$ contracted by the morphism $\phi_{i,j}$ satisfy $K_{X_{i,j-1}}\cdot C<0$; for any $i\in\{0,\dots,k-1\}$, the curves $C$ contracted by the morphism $\phi_{i}$ satisfy $K_{X_{i,j_i}}\cdot C<0$; and the the curves $C$ contracted by the morphism $\phi$ satisfy $K_{X_{k,j_k}}\cdot C<0$.
 \item For any $i\in\{0,\dots,k\}$ and $j\in\{1,\dots,j_i\}$, the curves $C$ contracted by the morphism $\phi_{i,j}^+$ satisfy $K_{X_{i,j}}\cdot C>0$.
 \item For any $i\in\{0,\dots,k-1\}$, the morphism $\phi_{i}$ contracts (at least) a $G$-stable divisor of $X_{i,j_i}$.
 \end{enumerate}
\end{prop}

\begin{proof}
\begin{enumerate} \item
Let $]a,b[$ or $[a,b[$ corresponding to an equivalent class in the family of polytopes $(Q^\epsilon)_{\epsilon\in\Qbb_{\geq 0}}$. In particular, there exists $I\subset\{1,\dots,m+r\}$ containing $K$ such that for any $\epsilon\in]a,b[$, $\tilde{Q}^\epsilon=P^\epsilon_I$ and $\epsilon\in\Omega_{K,I}^{max}$.

Denote by $X$ the variety given by the $G/H$-polytopes $Q^\epsilon$ with $\epsilon\in ]a,b[$, denote by $Y^b$ the variety given by the $G/H$-polytope $Q^b$ or to the $G/H'$-polytope $Q^b$ (as we define $Z$ in the previous section). And denote by $\phi$ the projective $G$-equivariant morphism from $X$ to $Y^b$.

Fix $c\in ]a,b[$.
Recall that a $G$-orbit of $X$ corresponding to a face $F_J^c$ of $Q^c$ with $J\subset I$, is sent to the $G$-orbit of $Y^b$ corresponding to the face $F_I^{b}$ of $Q^{b}$. 
Recall also that we describe the curves of horospherical varieties in Section~\ref{sec:curves}.
Hence, a curve $C_\mu$ of $X$ is contracted by $\phi$ if and only if the $G$-orbit of $X$ intersecting  $C_\mu$ in an open set (ie the $G$-orbit corresponding to the edge $\mu$ of $Q^{c}$) is sent to a closed $G$-orbit of $Y^b$; and a curve $C_{\alpha,v}$ of $X$ is contracted by $\phi$ if the closed $G$-orbit corresponding to $v$ is sent to a closed $G$-orbit isomorphic to $G/P'$ where $\alpha$ is a root of the parabolic subgroup $P'$. In other words, $C_\mu$ is contracted by $\phi$ if and only if for any $J\subset I$ such that $\mu=F_J^c$, the face $F_J^b$ of $Q^b$ is in fact a vertex of $Q^b$. And $C_{\alpha,v}$ is contracted by $\phi$ if and only if for any $J\subset I$ such that $v=F_J^c$ (which is not in $W_{\alpha,P}$), the vertex $F_J^b$ of $Q^b$ is in $W_{\alpha,P}$.

Let $D^c$ be the $\Qbb$-Cartier divisor defined by the moment polytope $Q^c$ and the pseudo-moment polytope $\tilde{Q}^c$. Then, for any $\epsilon\in[c,b[$ the $\Qbb$-Cartier (and ample) divisor defined by $(Q^\epsilon,\tilde{Q}^\epsilon)$ is $D^c+(\epsilon-c)K_X$.
But, by Proposition~\ref{prop:intersection}, $(D+(\epsilon-c)K_X)\cdot C_\mu$ is the integral length of the edge $\mu$ in $Q^\epsilon$ for any $\epsilon\in[c,b[$, we get by continuity that $C_\mu$ is contracted by $\phi$ if and only if $(D+(b-c)K_X)\cdot C_\mu=0$. Similarly, since $(D+(\epsilon-gamma)K_X)\cdot C_{\alpha,v}=\langle v,\alpha^\vee\rangle$ for any $\epsilon\in[c,b[$, we get that $C_{v,\alpha}$ is contracted by $\phi$ if and only if $(D+(b-c)K_X)\cdot C_{\alpha,v}=0$. In particular, for any curve $C$ contracted by $\phi$, we have $K_X\cdot C<0$.

\item By very similar arguments, we prove that if $]a,b[$ is an equivalent class in the family of polytopes  $(Q^\epsilon)_{\epsilon\in\Qbb_{\geq 0}}$, if $Y^a$ is the $G/H$-embedding associated to the $G/H$-polytope $Q^a$, and if $\phi^+:X\longrightarrow Y^a$ denote the $G$-equivariant morphism, then  for any curve $C$ contracted by $\phi^+$, we have $K_X\cdot C>0$.

\item Now consider the case where $]a,b[$ (or $[a,b[$) and $[b,b'[$ correspond to two successive equivalent classes in the family of polytopes $(Q^\epsilon)_{\epsilon\in\Qbb_{\geq 0}}$. The subset $I$ of $\{1,\dots,m+r\}$, the varieties $X$ and $Y^b$, and the morphism $\phi$ are defined as above. Here, by hypothesis, $b\not\in\Omega_{K,I}^{max}$ and there exists a proper subset $I'$ of $I$ containing $K$ such that $\tilde{Q}^b=P_{I'}^b$.
Then for any $c\in]a,b[$ and for any $i\in I\backslash I'\subset\{1,\dots,m\}$, $F_i^c$ is a facet of $Q^c$ (and corresponds to a $G$-stable divisor of $X$), but $F_i^b$ is not a facet of $Q^b$ (and corresponds to a $G$-stable divisor of $Y^b$ of codimension at least 2). Hence, $\phi$ contracts (at least) a $G$-stable divisor of $X$.
\end{enumerate}
 \end{proof}

\subsection{$\Qbb$-Gorenstein singularities}

In this section, we prove in particular that all the varieties $X_{i,j}$ defined in Section~\ref{sec:constructionvarietes} are $\Qbb$-Gorenstein.
We begin by giving a $\Qbb$-Gorenstein criterion in terms of moment polytopes.
\begin{prop}
Let $X$ be a projective $G/H$-embedding and let $D$ be an ample $\Qbb$-Cartier divisor. Denote by $\tilde{Q}$ the pseudo-moment polytope of $(X,D)$. Let $A$ and $\tilde{C}$ be the matrices defined in Section~\ref{sec:TheFamily}.
For any vertex $v$ of $\tilde{Q}$, we denote by $I_v$ the maximal subset of $\{1,\dots,m+r\}$ such that $v=F_{I_v}$.

Then $X$ is $\Qbb$-Gorenstein if and only if, for any vertex $v$ of $\tilde{Q}$, the linear system $A_{I_v}X=\tilde{C}_{I_v}$ have (at least) a solution.
\end{prop}

\begin{proof}
The proposition is just a translation, in terms of moment polytopes, of the criterion of $\Qbb$-Cartier divisor, applied to the divisor $K_X$.
\end{proof}

By applying this criterion to the family $(\tilde{Q}^\epsilon)_{\epsilon\in\Qbb_{\geq 0}}$, we easily get the following result.

\begin{cor}
Let $\epsilon\geq 0$ such that $Q^\epsilon$ is a $G/H$-polytope.
Let $X^\epsilon$ be the $G/H$-embedding corresponding to the $G/H$-polytope $Q^\epsilon$.
Denote by $I$ the subset of $I_0$ containing $K$, such that $\tilde{Q}^\epsilon=P_I^\epsilon$ and $\epsilon\in\Omega^{max}_{K,I}$. For any vertex $v^\epsilon$ of $\tilde{Q}$, we denote by $I_{v^\epsilon}$ the maximal subset of $I$ such that $v=F_{I_{v^\epsilon}}$.

Then $X^\epsilon$ is $\Qbb$-Gorenstein if and only if for any vertex $v^\epsilon$ of $\tilde{Q}^\epsilon$, the linear system $A_{I_{v^\epsilon}}X=\tilde{C}_{I_{v^\epsilon}}$ have (at least) a solution.
\end{cor}

Now, using Section~\ref{sec:poly}, we can know exactly when $X^\epsilon$ is $\Qbb$-Gorenstein (except for $Z$).

\begin{prop}
The varieties $X_{i,j}$ with $i\in\{0,\dots,k\}$ and $j\in\{0,\dots,j_i\}$ are $\Qbb$-Gorenstein.
And the varieties $Y_{i,j}$  with $i\in\{0,\dots,k\}$ and $j\in\{0,\dots,j_i\}$ are not $\Qbb$-Gorenstein.
\end{prop}

\begin{proof}
Let $i\in\{0,\dots,k\}$ and $j\in\{0,\dots,j_i\}$. The variety $X_{i,j}$ is defined by $G/H$-poltyopes $Q^\epsilon$ with $\epsilon$ in an open segment of $\Qbb_{\geq 0}$. In particular, for all these rational numbers $\epsilon$, for any vertex $v^\epsilon$ of $\tilde{Q}^\epsilon$, the linear system $A_{I_{v^\epsilon}}X=\tilde{B}_{I_{v^\epsilon}}+\epsilon\tilde{C}_{I_{v^\epsilon}}$ have a solution.
Hence, $A_{I_{v^\epsilon}}X=\tilde{C}_{I_{v^\epsilon}}$ has also a solution. It proves that $X_{i,j}$ is $\Qbb$-Gorenstein.

Now, let $i\in\{0,\dots,k\}$ and $j\in\{0,\dots,j_i\}$. The variety $Y_{i,j}$ is defined by the $G/H$-poltyope $Q^{\alpha_{i,j}}$. By Corollary~\ref{cor:closed}, there exist $J$ and $I$ subsets of $\{1,\dots,m+r\}$ such that $J\subset I$ and $\alpha_{i,j}$ is the extremity of the segment $\Omega_{J,I}^0$ or  $\{\alpha_{i,j}\}=\Omega_{J,I}^1$. But by Lemma~\ref{lem:separationextremities}, we can always choose $J$ and $I$ such that $\{\alpha_{i,j}\}=\Omega_{J,I}^1$. Moreover, by taking $J$ maximal with this property, we get a vertex $v$ of $Q^{\alpha_{i,j}}$ satisfying $\{\alpha_{i,j}\}=\Omega_{I_v,I}^1$. In particular, by Corollary~\ref{cor:closed}, the system $A_{I_v}X=\tilde{C}_{I_v}$ have no solution. It proves that $Y_{i,j}$ is not $\Qbb$-Gorenstein.
\end{proof}

\subsection{$\Qbb$-factorial singularities}

In this section, we prove that, for $D$ general, the MMP works in the family of projective $\Qbb$-factorial horospherical varieties.

\begin{prop}\label{prop:Qfactorial}
Suppose that $X$ is $\Qbb$-factorial.
Choose $D$ such that the vector $\tilde{B}$ is in the open set
$$\bigcup_{I\subset\{1,\dots,m+r\},\,|I|>n}\pi_I^{-1}(\Qbb^{|I|}\backslash\Im(A_I)),$$ where $\pi_I$ is the canonical projection of $\Qbb^{m+r}$ to its vector subspace of $\Qbb^{m+r}$ corresponding to the coordinates in $I$.

Then, for any $i\in\{0,\dots,k\}$ and $j\in\{0,\dots,j_i\}$, the variety $X_{i,j}$ is $\Qbb$-factorial. 
\end{prop}

\begin{proof}
A horospherical variety $X$ is $\Qbb$-factorial if and only if the colored cones of $\mathbb{F}_X$ are all simplicial and any color of $X$ is the unique color of a colored edge of $\mathbb{F}_X$. In terms of moment polytopes, if $Q$ is a moment polytope of $X$, then $X$ is $\Qbb$-factorial if and only if $Q$ is simple (ie each vertex belongs exactly to $n$ facets, where $n$ is the dimension of $Q$), $Q$ intersects a wall $W_{\alpha,P}$, with $\alpha\in S\backslash R$, only along one of its facets and a facet is never in 2 walls $W_{\alpha,P}$ and $W_{\beta,P}$ with $\alpha\neq\beta\in S\backslash R$.

Let $i\in\{0,\dots,k-1\}$ and $j\in\{0,\dots,j_i\}$. Let $\epsilon\in]\alpha_{i,j},\alpha_{i,j+1}[$ so that $Q^\epsilon$ is a moment polytope of $X_{i,j}$. Let $I\subset I_0$ containing $K$ such that $\tilde{Q}^\epsilon=P_I^\epsilon$ and $\epsilon\in\Omega^{max}_{K,I}$. 

Let $v$ be a vertex of $Q^\epsilon$. Denote by $I_v$ the maximal subset of $I$ such that $v=F_{I_v}^\epsilon$. Note that $|I_v|\geq n$. We want to prove that $|I_v|=n$.
In particular, $\epsilon$ is a point of $\Omega^1_{I_v}$, which is open (because contains $]\alpha_{i,j},\alpha_{i,j+1}[$). It implies that, for any $\eta\in\Qbb$,  $\tilde{B}_{I_v}+\eta\tilde{C}_{I_v}$ is in the image of $A_{I_v}$. In particular, $\tilde{B}_{I_v}$ is in the image of $A_{I_v}$. By hypothesis on $D$, the cardinality of $I_v$ has to be $n$. This proves that $X_{i,j}$ is $\Qbb$-factorial.
\end{proof}

\begin{rem}
The open set where $\tilde{B}$ is chosen, is clearly not empty and dense in $\Qbb^{m+r}$, because for any $I$ of cardinality greater than $n$, the image of $A_I$ is of codimension at least one. But, since $X$ is $\Qbb$-factorial, any vector $\tilde{B}\in\Qbb^{m+r}$ gives a $\Qbb$-Cartier divisor.
\end{rem}

Taking $D$ general, we also get the following result.

\begin{prop}\label{prop:onlyrays}
Suppose that $X$ is $\Qbb$-factorial.
If $D$ is general in the set of ample $\Qbb$-Cartier divisors, all morphisms $\phi_{i,j}$, $\phi_{i,j}^+$, $\phi_i$ and $\phi$ defined in Section~\ref{sec:constructionvarietes} are contractions of rays of the corresponding effective cones $NE(X_{i,j})$. 
\end{prop}

\begin{proof}
Let $Y$ be a projective $\Qbb$-factorial $G/H$-embedding whose colored fan is made from some edges generated by $x_i$ with $i\in\{1,\dots,m\}$ and $\alpha_M^\vee$ with $\alpha\in S\backslash R$. In particular, the vector space $\operatorname{Cartier}(X)_\Qbb$ of $\Qbb$-Cartier divisors of $X$ projects on the vector space $\operatorname{Cartier}(Y)_\Qbb$ of $\Qbb$-Cartier divisors of $Y$. Denote by $p_Y$ this projection.
Let $V$ be a vector subspace of $N_1(Y)_\Qbb$ of dimension at least~2. Denote by $V^*$ the dual of $V$ in $N^1(X)_\Qbb$. Then $V^*+\Qbb K_Y$ is a proper vector subspace of $N^1(X)_\Qbb$. Hence, the intersection $\mathcal{I}_Y$, for all faces of $NE(X)$ of dimension at least~2 of direction the vector subspace $V$, of $N^1(Y)\backslash (V^*+\Qbb K_Y)$ is a dense open set of $N^1(Y)_\Qbb$. If we denote by $q_Y$ the projection of $\operatorname{Cartier}(Y)_\Qbb$ on $N^1(Y)_\Qbb$, then $(q_Y\circ p_Y)^{-1}(\mathcal{I}_Y)$ is open and dense in $\operatorname{Cartier}(X)_\Qbb$.

Consider now, an ample $\Qbb$-Cartier divisor $D$ in the dense open intersection of the set given in Proposition~\ref{prop:Qfactorial} with the sets $(q_Y\circ p_Y)^{-1}(\mathcal{I}_Y)$ for all projective $\Qbb$-factorial $G/H$-embeddings $Y$ whose colored fan is made from some edges generated by $x_i$ with $i\in\{1,\dots,m\}$ and $\alpha_M^\vee$ with $\alpha\in S\backslash R$.
Then, for any $i\in\{0,\dots,k\}$ and $j\in\{1,\dots,j_i\}$, the divisor $p_{X_{i,j-1}}(D)+\alpha_{i,j}K_{X_{i,j-1}}$ vanishes at most on a ray of $NE(X_{i,j-1})$, so that  $\phi_{i,j}$ is the contraction of a ray of $NE(X_{i,j-1})$. Similarly, for any $i\in\{0,\dots,k\}$ and $j\in\{1,\dots,j_i\}$, the divisor $p_{X_{i,j}}(D)+\alpha_{i,j}K_{X_{i,j}}$ vanishes at most on a ray of $NE(X_{i,j})$, so that  $\phi_{i,j}^+$ is the contraction of a ray of $NE(X_{i,j})$. And for any $i\in\{0,\dots,k\}$, the divisor $p_{X_{i,j_i}}(D)+\alpha_{i,j_i+1}K_{X_{i,j_i}}$ vanishes at most on a ray of $NE(X_{i,j_i})$, so that  $\phi_i$ ($\phi$ if $i=k$) is the contraction of a ray of $NE(X_{i,j_i})$.
\end{proof}

This result seems to be true without the hypothesis of $\Qbb$-factoriality, nevertheless the proof above really uses this hypothesis.

\subsection{General fibers of Mori fibrations}\label{sec:generalfiber}

In this section, we study the general fibers of the morphism $\phi:X_{k,j_{k}}\longrightarrow Z$ defined in Section~\ref{sec:constructionvarietes} and prove the last statement of Theorem~\ref{th:maingeneral}. To simplify the notations, we suppose that $k=1$ and $j_0=0$, in particular $X_{k-1,j_{k-1}}=X$. We also denote by $\epsilon_1$ the rational number $\alpha_{1,0}$.
We distinguish two cases.\\

First suppose that $Q^{\epsilon_1}$ and $Q$ have the same dimension. It implies that $Q^{\epsilon_1}$ is in some wall $W_{\alpha,P}$ with $\alpha\in S\backslash R$, ie $R^1\neq R$. 
By Proposition~\ref{prop:onlyrays}, if $X$ is $\Qbb$-factorial, then for $D$ general, we have $|R^1\backslash R|=1$ and then $P^1/P$ of Picard number~1. Indeed, if $|R^1\backslash R|>1$, $\phi$ is not a contraction of an extremal ray, because $G/P\longrightarrow G/P^1$ is clearly not a contraction of an extremal ray.
In that case, the fibers of the morphism $\phi$ is a flag variety of Picard number one.\\

Secondly, suppose that the dimension of $Q^{\epsilon_1}$ is less than the dimension of $Q$.

Let $x_0$ (resp. $x_0^1$) be the unique point of the open $G$-orbit of $X$ (resp. $Z$) fixed by $H$ (resp. $H^1$). We claim that the general fiber of $\phi$ is the closure in $X$ of the $H^1$-orbit $H^1\cdot x_0$. Indeed, since $G\cdot x_0^1$ is open in $Z$, the fiber $\phi^{-1}(x_0^1)$ is a general fiber of $\phi$. This fiber intersected by the open $G$-orbit $G\cdot x_0$ of $X$ is $H^1\cdot x_0$. We conclude the claim by the fact that, since $X$ is irreducible, a general fiber of $\phi$ is also irreducible. Denote by $F_\phi$ the $H^1$-variety $\phi^{-1}(x_0^1)$.

Since an open set of $X$ is isomorphic to the bundle $G\times^{H^1}F_\phi$, we get that $F_\phi$ has the same singularities as $X$ (normal and $\Qbb$-Gorenstein, and $\Qbb$-factorial or smooth according to the hypothesis done on $X$).
Moreover, the unipotent radical $R_u(H^1)$ of $H^1$ (which is also the unipotent radical of $P^1$) acts trivially on $x_0$, then it also acts trivially on $H^1\cdot x_0$ and $F_phi$. Hence, $F_\phi$ is a $L^1$-variety, where $L^1=H^1/R_u(H^1)$ is reductive. Let $H^2:=H/R_u(H^1)$. Then the open $L^1$-orbit of $F_\phi$ is isomorphic to $L^1/H^2$ and is horospherical, because $H^2$ contains the unipotent radical $U/R_u(H^1)$ of the Borel subgroup $B/R_u(H^1)$ of $L^1$. In the rest of the section, we describe the $L^1/H^2$-embedding $F_\phi$.

First describe the combinatorial data associated to the the horospherical homogeneous space $L^1/H^2$.
The simple roots of $L^1$ are those of $P^1$, ie the simple roots in the set $R^1$. Then the normalizer $P^2$ of $H^2$ in $L^1$ is the parabolic subgroup of $L^1$ whose set of simple roots is $R$. Moreover, the set of characters of $P^2$ trivial on $H^2$ is isomorphic to the quotient of the set of characters of $P$ trivial on $H$ with $M^1$. Then we set $M^2:=M/M^1$. The set of colors of the horospherical homogeneous space $L^1/H^2$ is $R^1\backslash R$.  

Let $I$ be the maximal subset of $\{1,\dots,m+r\}$ such that $\tilde{Q}^{\epsilon_1}\subset\{x\in M\,\mid\, A_{I}x=\tilde{B}_I+\epsilon_1\tilde{C}_I\}$. 
In particular, $M^1_\Qbb$ equals the kernel of $A_{I}$. 
Now, we prove that the projection of $\tilde{Q}$ on $M^2_\Qbb:=M_\Qbb/M^1_\Qbb$  is the polytope $\tilde{Q}^2:=\{x\in M_\Qbb/\Ker(A_{I})\,\mid\, A_{I}x\geq \tilde{B}_I\}$. 
One inclusion is obvious. To prove the second one, let $x+\Ker(A_{I})\in\tilde{Q}^2$ be a vertex $F^2_{J}$ of $\tilde{Q}^2$ where $J\subset I$. 
Then $A_{J}x=\tilde{B}_{J}$ and $A_{I}x\geq\tilde{B}_{I}$. 
Moreover, by maximality of $I$, there exists $x'\in M_\Qbb$ such that $A_{I}x'=\tilde{B}_{I}+\epsilon_1\tilde{C}_{I}$ and $A_{\bar{I}}x'>\tilde{B}_{\bar{I}}+\epsilon_1\tilde{C}_{\bar{I}}$ (where $\bar{I}=\{1,\dots,m+r\}\backslash I$). 
Hence, for $\lambda>0$ small enough, we have $x'':=\lambda x+(1-\lambda)x'$ satisfies $A_{J}x''=\tilde{B}_{J}+\epsilon''\tilde{C}_{J}$, $A_{I}x''\geq\tilde{B}_{I}+\epsilon''\tilde{C}_{I}$ and $A_{\bar{I}}x''>\tilde{B}_{\bar{I}}+\epsilon''\tilde{C}_{\bar{I}}$, with $\epsilon''=(1-\lambda)\epsilon_1$, so that $x''$ is a point of the face $F^{\epsilon''}_{J}$ of $\tilde{Q}^{\epsilon''}$. 
Since $0\leq\epsilon''<\epsilon_1$, the face $F_{J}$ of $\tilde{Q}$ is not empty. 
But every point of $F_{J}$ projects on $F^2_{J}$, which is a point. Then, every vertex of $\tilde{Q}^2$ is the image of a point of $\tilde{Q}$ by the projection on $M^2_\Qbb$, that proves the second inclusion.

By translation, we get that the projection of $Q$ on $X(P)_\Qbb/M^1_\Qbb$ is the polytope $Q^2:=\{x\in \bar{v^0}+M_\Qbb/\Ker(A_{I})\,\mid\, A_{I}x\geq B_I\}$, where $\bar{v^0}$ is the image of $v^0$ in $X(P)_\Qbb/M^1_\Qbb$.

Suppose now that $D$ is Cartier and very ample (or replace $D$ by a multiple of $D$, see Remark~\ref{rem:veryample}). Then, by Proposition~\ref{prop:veryample}, $X$ is isomorphic to the closure of $G\cdot [\sum_{\chi\in (v^0+M)\cap Q}v_\chi]$ in $\Pbb(\oplus_{\chi\in (v^0+M)\cap Q}V(\chi))$. Then $F_\phi$ is isomorphic to $L^1\cdot [\sum_{\chi\in (v^0+M)\cap Q}v_\chi]$ in $\Pbb(\oplus_{\chi\in (v^0+M)\cap Q}V_{L^1}(\chi))$. But, for any $\chi\in X(P)$ and any $\chi'\in M^1$, the $L^1$-modules $V_{L^1}(\chi)$ and $V_{L^1}(\chi+\chi')$ are isomorphic, so that $Y$ is isomorphic to $L^1\cdot [\sum_{\chi\in (\bar{v^0}+M^2)\cap Q^2}v_\chi]$ in $\Pbb(\oplus_{\chi\in (\bar{v^0}+M^2)\cap Q^2}V_{L^1}(\chi))$. It proves that $F_\phi$ is the $L^1/H^2$-embedding associated to the polytope $Q^2$. Remark that $Q^2$ is of maximal dimension in $M^2_\Qbb$ because $Q$ is of maximal dimension in $M_\Qbb$.

Suppose now that $X$ is $\Qbb$-factorial. We prove that, for $D$ general, $I$ is of cardinality $\codim(Q^{\epsilon_1})+1=\dim(Q^2)+1$, and that $Q^2$ is a simplex of $M^2_\Qbb$ such that all wall $W^2_\alpha$ in $X(P^2)_\Qbb$ with $\alpha\in R^1\backslash R$ gives a facet of $Q^2$.
Suppose that $D$ (and then $\tilde{B}$) satisfies the following condition: for any subset $J$ of $\{1,\dots,m+r\}$ with $|J|\geq\dim(\Im(A_J))+2$, $\tilde{B}_{J}$ is not in the (proper) vector subspace of $\Qbb^{m+r}$ generated by $\Im(A_J)$ and $\tilde{C}_J$. Such a $D$ is general, since $X$ is $\Qbb$-factorial.
Note that $\codim(Q^{\epsilon_1})=\codim(\Ker(A_{I}))=\dim(\Im(A_{I}))$. Now, if $|I|>\codim(Q^{\epsilon_1})+1$ then by hypothesis on $D$, $\tilde{B}_{I}+\epsilon \tilde{C}_{I}$ is never in $\Im(A_{I})$, which contradicts that $F_{I}^{\epsilon_1}=Q^{\epsilon_1}$ is not empty.
Moreover,  $Q^2$ is a polytope (and then bounded) of maximal dimension in $M^2_\Qbb$, so the number $|I|$ of inequations defining $Q^2$ must be at least $\dim(Q^2)+1$.
Hence, $|I|=\dim(Q^2)+1$ and then $Q^2$ is a simplex of $M^2_\Qbb$.
The set $R^1\backslash R$ of colors of $L^1/H^2$ is clearly contained in $I$ by definition of $I$ and $R^1$, so that, for any $\alpha\in R^1\backslash R$, $W^2_\alpha\cap Q^2$ is a face of $Q^2$. But, since each $i\in I$ correspond necessarily to a facet of $Q^2$ because $|I|=\dim(Q^2)+1$, this latter face is necessarily a facet.

In that case, $F_\phi$ is a projective $\Qbb$-factorial $L^1/H^2$-embedding of Picard number~1 (cf \cite[Eq. (4.5.1)]{these} for an explicit formula of the Picard number of $\Qbb$-factorial horospherical varieties).

\section{Examples}\label{sec:examples}

In this section we give 5 examples with $G=\operatorname{SL}_3$. The first three ones give the MMP for the same horospherical smooth variety but with a different ample Cartier divisor. The forth one gives a flip consisting of exchanging colors (see \cite[Section 4.5]{brionmori}). And the last one give the MMP for a $\Qbb$-Gorenstein (not $\Qbb$-factorial) variety, we observe in particular a flip from a divisorial contraction.

Fix a Borel subgroup $B$ of $G$. Denote by $\alpha$ and $\beta$ the two simple roots of $G$.

\begin{ex}\label{ex:horo1}
 Consider the horospherical subgroup $H$ defined as the kernel of the character $\varpi_\alpha+\varpi_\beta$ of $B$. In that case we have $N$ and $M$ isomorphic to $\Zbb$. The horospherical homogeneous space has two colors $\alpha$ and $\beta$ whose image in $N$ are respectively $\alpha_M^\vee=\beta_M^\vee=1$.

If $\chi$ is a character of $B$, we denote by $\Cbb_\chi$ the line $\Cbb$ where $B$ acts by $b\cdot z=\chi(b)z$ for any $b\in B$ and $z\in\Cbb$.

Let $X$ be the $\Pbb^1$-bundle $X=G\times^B\Pbb(\Cbb_0\oplus\Cbb_{\varpi_\alpha+\varpi_\beta})$ over $G/B$, it is a smooth $G/H$-embedding. Its colored fan is the unique complete fan (of dimension~1) without color. Denote by $X_1$ and $X_2$ the two irreducible $G$-stable divisors of $X$, respectively corresponding to the primitive elements $x_1=1$ and $x_2=-1$ of $N$. Here $-K_X=X_1+X_2+2D_\alpha+2D_\beta$.

Choose $D=X_1+2X_2+2D_\alpha+2D_\beta$. Then the moment polytope $Q$ is the segment $[\varpi_\alpha+\varpi_\beta, 4(\varpi_\alpha+\varpi_\beta)]$ in the dominant chamber of $(G,B)$.

The family $(Q^\epsilon)_{\epsilon\geq 0}$ is given by:
\begin{itemize}
\item for any $\epsilon\in[0,1[$, $Q^\epsilon$ is the segment $[(1-\epsilon)(\varpi_\alpha+\varpi_\beta), (4-3\epsilon)(\varpi_\alpha+\varpi_\beta)]$;
\item for any $\epsilon\in[1,\frac{4}{3}[$, $Q^\epsilon$ is the segment $[0,(4-3\epsilon)(\varpi_\alpha+\varpi_\beta)]$;
\item $Q^\frac{4}{3}$ is the point $0$.
\end{itemize}

Hence, the MMP from $(X,D)$ gives a divisorial contraction from $X$ to the projective $G/H$-embedding with the two colors $\alpha$ and $\beta$, which is not $\Qbb$-factorial but $\Qbb$-Gorenstein, and Fano. It finishes by a Mori fibration from this Fano variety to a point.

Note that, here, the divisorial contraction contracts 2 divisors: the zero and infinite section of the $\Pbb^1$-bundle $X$. In particular, it is not the contraction of a ray of $NE(X)$.
\end{ex}

\begin{ex}\label{ex:horo2}
We keep the same $G/H$-embedding $X$ but we choose another ample divisor $D=X_1+2X_2+3D_\alpha+2D_\beta$. Then the moment polytope $Q$ is the segment $[2\varpi_\alpha+\varpi_\beta, 5\varpi_\alpha+4\varpi_\beta]$.

The family $(Q^\epsilon)_{\epsilon\geq 0}$ is given by:
\begin{itemize}
\item for any $\epsilon\in[0,1[$, $Q^\epsilon$ is the segment $[2\varpi_\alpha+\varpi_\beta-\epsilon(\varpi_\alpha+\varpi_\beta), 5\varpi_\alpha+4\varpi_\beta-3\epsilon(\varpi_\alpha+\varpi_\beta)]$;
\item for any $\epsilon\in[1,\frac{4}{3}[$, $Q^\epsilon$ is the segment $[\varpi_\alpha,5\varpi_\alpha+4\varpi_\beta-3\epsilon(\varpi_\alpha+\varpi_\beta)]$;
\item $Q^\frac{4}{3}$ is the point $\varpi_\alpha$.
\end{itemize}

Hence, the MMP from $(X,D)$ gives a divisorial contraction from $X$ to the projective $G/H$-embedding with the color $\beta$. This is a contraction of the ray of $NE(X)$ generated by $C_{\beta,X_1}$.
It finishes by a Mori fibration from this $\Qbb$-factorial variety to the flag variety $G/P_\alpha$.
\end{ex}

\begin{ex}\label{ex:horo3}
We still keep the same $G/H$-embedding $X$ and we choose now the ample divisor $D=X_2+2D_\alpha+2D_\beta$. Then the moment polytope $Q$ is the segment $[2(\varpi_\alpha+\varpi_\beta), 3(\varpi_\alpha+\varpi_\beta)]$.

The family $(Q^\epsilon)_{\epsilon\geq 0}$ is given by:
\begin{itemize}
\item for any $\epsilon\in[0,\frac{1}{2}[$, $Q^\epsilon$ is the segment $[2(\varpi_\alpha+\varpi_\beta)+\epsilon(-\varpi_\alpha+\varpi_\beta),3(\varpi_\alpha+\varpi_\beta)-3\epsilon(\varpi_\alpha+\varpi_\beta)]$;
\item $Q^\frac{1}{2}$ is the point $\frac{3}{2}(\varpi_\alpha+\varpi_\beta)$.
\end{itemize}

Hence, the MMP from $(X,D)$ gives the Mori fibration $X\longrightarrow G/B$, whose fibers are projective lines. 
\end{ex}

In the following example, we illustrate a flip consisting of exchanging colors.

\begin{ex}\label{ex:horo4}
Consider the horospherical subgroup $H$ defined as the kernel of the character $\varpi_\alpha+2\varpi_\beta$ of $B$. In that case we have $N$ and $M$ isomorphic to $\Zbb$. The horospherical homogeneous space has two colors $\alpha$ and $\beta$ whose image in $N$ are respectively $\alpha_M^\vee=1$ and $\beta_M^\vee=2$.

Let $X$ be the $\Qbb$-factorial $G/H$-embedding whose colored fan is the complete fan with color $\beta$. Denote by $X_1$ the irreducible $G$-stable divisor of $X$, corresponding to the primitive elements $x_1=-1$. Here $-K_X=X_1+2D_\alpha+2D_\beta$.

Consider $D=3X_1+2D_\alpha+2D_\beta$. Then the moment polytope $Q$ is the segment $[\varpi_\alpha, 5\varpi_\alpha+8\varpi_\beta]$.

The family $(Q^\epsilon)_{\epsilon\geq 0}$ is given by:
\begin{itemize}
\item for any $\epsilon\in[0,1[$, $Q^\epsilon$ is the segment $[(1-\epsilon)\varpi_\alpha,5\varpi_\alpha+8\varpi_\beta-\epsilon(3\varpi_\alpha+4\varpi_\beta)]$;
\item $Q^1$ is the segment $[0,2\varpi_\alpha+4\varpi_\beta]$;
\item for any $\epsilon\in]1,\frac{5}{3}[$, $Q^\epsilon$ is the segment $[2(\epsilon-1)\varpi_\beta,5\varpi_\alpha+8\varpi_\beta-\epsilon(3\varpi_\alpha+4\varpi_\beta)]$;
\item $Q^\frac{5}{3}$ is the point $\frac{4}{3}\varpi_\beta$.
\end{itemize}

Hence, the MMP from $(X,D)$ first gives a flip $X\longrightarrow Y\longleftarrow X^+$, where $Y$ is the $G/H$-embedding corresponding to the complete colored fan with the two colors $\alpha$ and $\beta$ and $X^+$ is the $G/H$-embedding corresponding to the complete colored fan with the color $\alpha$.
It finishes by a Mori fibration from $X^+$ to the flag variety $G/P_\beta$.
\end{ex}

\begin{ex}\label{ex:horo5}
Consider the case where the horospherical subgroup $H$ is the maximal unipotent $U$ subgroup of $B$. Then the lattice $M$ is the lattice of characters of $B$ with basis $(\varpi_\alpha,\varpi_\beta)$, and $N$ is the coroot lattice with basis $(\alpha^\vee,\beta^\vee)$. Here $\alpha^\vee=\alpha^\vee_M$ and $\beta^\vee=\beta^\vee_M$.

Let $X$ be the $G/H$-embedding whose colored fan is the complete colored fan with color $\beta$, and edges generated by $x_1:=-\beta^\vee$, $x_2:=\alpha^\vee$ and $x_3:=\beta^\vee-\alpha^\vee$. Note that $X$ is not $\Qbb$-factorial because $\beta^\vee$ is not in an edge of $\mathbb{F}_X$.  An anticanonical divisor of $X$ is $X_1+X_2+X_3+2D_\alpha+2D_\beta$, and it is Cartier so that $X$ is $\Qbb$-Gorenstein.

Consider $D=3X_1+X_3+D_\alpha+D_\beta$. Then $D$ is an ample Cartier divisor whose moment polytope $Q$ is the triangle with vertices $\varpi_\alpha$, $\varpi_\alpha+4\varpi_\beta$ and $5\varpi_\alpha+4\varpi_\beta$. The matrices defining the family $(Q^\epsilon)_{\epsilon\geq 0}$ are $$A=\left(\begin{array}{cc}
0 & -1\\
1 & 0\\
-1 & 1\\
1 & 0\\
0 & 1
\end{array}\right),\,B=\left(\begin{array}{c}
-4\\
1\\
-1\\
0\\
0
\end{array}\right)\mbox{ and }C=\left(\begin{array}{c}
3\\
-1\\
1\\
0\\
0
\end{array}\right).$$

Then the family $(Q^\epsilon)_{\epsilon\geq 0}$ is given by:
\begin{itemize}
\item for any $\epsilon\in[0,1[$, $Q^\epsilon$ is the triangle with vertices $(1-\epsilon)\varpi_\alpha$, $\varpi_\alpha+4\varpi_\beta -\epsilon(\varpi_\alpha+3\varpi_\beta)$ and $5\varpi_\alpha+4\varpi_\beta-\epsilon(4\varpi_\alpha+3\varpi_\beta)$; 
\item $Q^1$ is the triangle with vertices $0$, $\varpi_\beta$ and $\varpi_\alpha+\varpi_\beta$;
\item for any $\epsilon\in]1,\frac{5}{4}[$, $Q^\epsilon$ is the triangle with vertices $(\epsilon-1)\varpi_\beta$, $(4-3\epsilon)\varpi_\beta$ and $5\varpi_\alpha+4\varpi_\beta-\epsilon(4\varpi_\alpha+3\varpi_\beta)$;
\item $Q^\frac{5}{4}$ is the point $\frac{1}{4}\varpi_\beta$.
\end{itemize}
\begin{figure}[htbp]
\begin{center}
\input{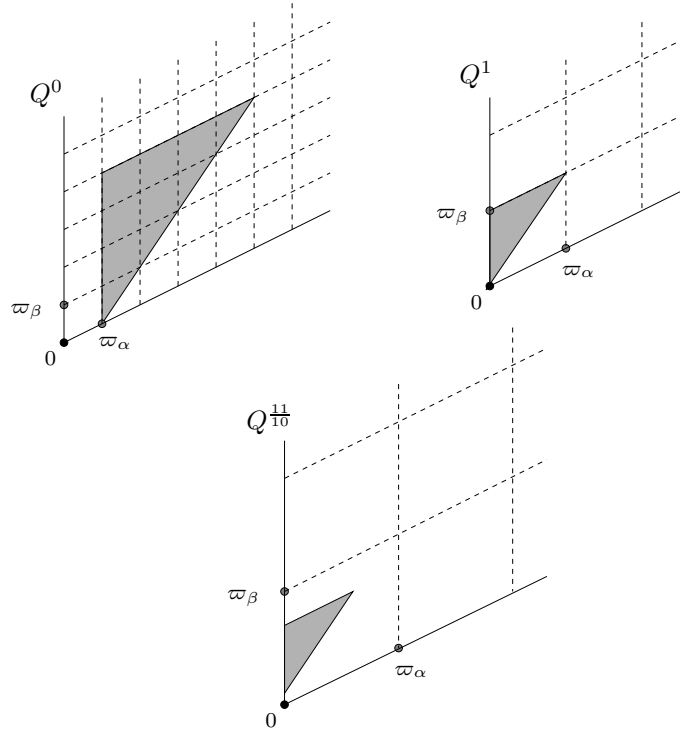}
\caption{Evolution of $Q^\epsilon$ in Example \ref{ex:horo5}}
\label{figure:exemple7}
\end{center}
\end{figure}
We illustrate the first three classes of the family $(Q^\epsilon)_{\epsilon\geq 0}$ in Figure~\ref{figure:exemple7}.

Hence, the MMP from $(X,D)$ first gives a flip $X\longrightarrow Y\longleftarrow X^+$, where $Y$ is the $G/H$-embedding corresponding to the same complete colored fan as $X$ but with the two colors $\alpha$ and $\beta$ and $X^+$ is the $G/H$-embedding corresponding to the same complete colored fan as $X$ but with the color $\alpha$ (instead of $\beta$). Note that the map $X\longrightarrow Y$ contracts the divisor $X_2$ and that the map $X^+\longrightarrow Y$ does not contract a divisor.
It finishes by a Mori fibration from $X^+$ to the flag variety $G/P_\beta$.

\end{ex}

\bibliographystyle{amsalpha}
\bibliography{MMPviapoly}

\providecommand{\bysame}{\leavevmode\hbox to3em{\hrulefill}\thinspace}
\providecommand{\MR}{\relax\ifhmode\unskip\space\fi MR }
\providecommand{\MRhref}[2]{%
  \href{http://www.ams.org/mathscinet-getitem?mr=#1}{#2}
}
\providecommand{\href}[2]{#2}
\begin{thebibliography}{Oda88}

\bibitem[Bri89]{briondiv}
Michel Brion, \emph{Groupe de picard et nombres caract\'eristiques des
  vari\'et\'es sph\'eriques}, Duke Math. J. \textbf{58} (1989), no.~2,
  397--424.

\bibitem[Bri93]{brionmori}
\bysame, \emph{Vari\'et\'es sph\'eriques et th\'eorie de {M}ori}, Duke Math. J.
  \textbf{72} (1993), no.~2, 369--404.

\bibitem[Ful93]{fulton}
William Fulton, \emph{Introduction to toric varieties}, Annals of Mathematics
  Studies, vol. 131, Princeton University Press, Princeton, NJ, 1993, The
  William H. Roever Lectures in Geometry.

\bibitem[Kno91]{knop}
Friedrich Knop, \emph{The {L}una-{V}ust theory of spherical embeddings},
  Proceedings of the {H}yderabad {C}onference on {A}lgebraic {G}roups
  ({H}yderabad, 1989) (Madras), Manoj Prakashan, 1991, pp.~225--249.

\bibitem[Mat02]{matsuki}
Kenji Matsuki, \emph{Introduction to the {M}ori program}, Universitext,
  Springer-Verlag, New York, 2002.

\bibitem[Oda88]{oda}
Tadao Oda, \emph{Convex bodies and algebraic geometry}, Ergebnisse der
  Mathematik und ihrer Grenzgebiete (3) [Results in Mathematics and Related
  Areas (3)], vol.~15, Springer-Verlag, Berlin, 1988, An introduction to the
  theory of toric varieties, Translated from the Japanese.

\bibitem[Pas06]{these}
Boris Pasquier, \emph{Vari\'et\'es horosph\'eriques de {F}ano}, Ph.D. thesis,
  Universit\'e Joseph Fourier, Grenoble~1, Available at
  http://tel.archives-ouvertes.fr/tel-00111912, 2006.

\bibitem[Pas08]{Fanohoro}
\bysame, \emph{Vari\'et\'es horosph\'eriques de {F}ano}, Bull. Soc. Math.
  France \textbf{136} (2008), no.~2, 195--225.

\bibitem[Pas09]{2orbits}
\bysame, \emph{On some smooth projective two-orbit varieties with {P}icard
  number 1}, Math. Ann. \textbf{344} (2009), no.~4, 963--987.

\bibitem[Rei83]{Reid}
Miles Reid, \emph{Decomposition of toric morphisms}, Arithmetic and geometry,
  {V}ol. {II}, Progr. Math., vol.~36, Birkh\"auser Boston, Boston, MA, 1983,
  pp.~395--418.

\end{thebibliography}

\end{document}